\documentclass[10pt,letterpaper]{amsart}
\usepackage{charter}
\usepackage[OT2,OT1]{fontenc}
\usepackage[latin9]{luainputenc}
\usepackage{array}
\usepackage{float}
\usepackage{mathrsfs}
\usepackage{bm}
\usepackage{multirow}
\usepackage{amstext}
\usepackage{amsthm}
\usepackage{amssymb}
\usepackage{cancel}
\usepackage{stmaryrd}

\makeatletter

\@ifundefined{pageheight}{\let\pageheight\pdfpageheight}{}
\@ifundefined{pagewidth}{\let\pagewidth\pdfpagewidth}{}
\pageheight\paperheight
\pagewidth\paperwidth

\providecommand{\tabularnewline}{\\}

\numberwithin{equation}{section}
\numberwithin{figure}{section}
\theoremstyle{plain}
\newtheorem{thm}{\protect\theoremname}[section]
\theoremstyle{definition}
\newtheorem{defn}[thm]{\protect\definitionname}
\theoremstyle{remark}
\newtheorem{rem}[thm]{\protect\remarkname}
\theoremstyle{plain}
\newtheorem{prop}[thm]{\protect\propositionname}
\theoremstyle{plain}
\newtheorem{cor}[thm]{\protect\corollaryname}
\theoremstyle{plain}
\newtheorem{lem}[thm]{\protect\lemmaname}
\theoremstyle{definition}
\newtheorem{example}[thm]{\protect\examplename}

\usepackage{amsfonts}
\usepackage{mathrsfs}
\usepackage{graphicx}
\usepackage{amscd}
\usepackage{url}

\newcommand{\cyr}{
\renewcommand\rmdefault{wncyr} \renewcommand\sfdefault{wncyss} \renewcommand\encodingdefault{OT2} \normalfont
\selectfont
}
\DeclareTextFontCommand{\textcyr}{\cyr}    

   \def\@settitle
{
\begin{center} \baselineskip14\p@\relax \LARGE\@title \end{center}
}

\usepackage[normalem]{ulem}
\usepackage[mathscr]{eucal}
\numberwithin{equation}{section}

\usepackage{verbatim}
\usepackage{amscd}
\usepackage{bm}\usepackage[mathscr]{eucal}
\usepackage[all]{xy}
\usepackage{enumerate}
\usepackage{color}
\usepackage{colortbl}

\input xyimport.tex

\title{\bf Key varieties for prime $\mathbb{Q}$-Fano threefolds \\ defined by Freudenthal triple systems}

\author{Hiromichi Takagi}

\address{Department of Mathematics, Gakushuin University, 
Mejiro, Toshima-ku, Tokyo 171-8588, Japan}
\email{hiromici@math.gakushuin.ac.jp}

\newcommand{\sO}{\mathcal{O}}

\newcommand{\mA}{\mathbb{A}}
\newcommand{\mC}{\mathbb{C}}

\newcommand{\mP}{\mathbb{P}}
\newcommand{\mQ}{\mathbb{Q}}

\newcommand{\rank}{\mathrm{rank}\,}

\numberwithin{equation}{section}

 \newcounter{myparagraph}[subsection]

\makeatother

\providecommand{\corollaryname}{Corollary}
\providecommand{\definitionname}{Definition}
\providecommand{\examplename}{Example}
\providecommand{\lemmaname}{Lemma}
\providecommand{\propositionname}{Proposition}
\providecommand{\remarkname}{Remark}
\providecommand{\theoremname}{Theorem}

\begin{document}
\maketitle 
\begin{abstract}
In this paper, we are concerned with the classification of complex
prime $\mQ$-Fano $3$-folds of anti-canonical codimension 4 which
are produced, as weighted complete intersections of appropriate weighted
projectivizations of certain affine varieties related with $\mP^{1}\times\mP^{1}\times\mP^{1}$-fibrations.
Such affine varieties or their appropriate weighted projectivizations
are called key varieties for prime $\mQ$-Fano 3-folds. We realize
that the equations of the key varieties can be described conceptually
by Freudenthal triple systems (FTS, for short). The paper consists
of two parts. In Part 1, we revisit the general theory of FTS; the
main purpose of Part 1 is to derive the conditions of so called strictly
regular elements in FTS so as to fit with our description of key varieties.
Then, in Part 2, we define several key varieties for prime $\mQ$-Fano
3-folds from the conditions of strictly regular elements in FTS. Among
other things obtained in Part 2, we show that there exists a $14$-dimensional
factorial affine variety $\mathfrak{U}_{\mA}^{14}$ of codimension
$4$ in an affine $18$-space with only Gorenstein terminal singularities,
and we construct examples of prime $\mQ$-Fano $3$-folds of No.20544
in \cite{GRDB} as weighted complete intersections of the weighted
projectivization of $\mathfrak{U}_{\mA}^{14}$ in the weighted projective
space $\mP(1^{15},2^{2},3)$. We also clarify in Part 2 a relation
between $\mathfrak{U}_{\mA}^{14}$ and the $G_{2}^{(4)}$-cluster
variety, which is a key variety for prime $\mQ$-Fano 3-folds constructed
in \cite{CD}. 
\end{abstract}

\maketitle
\markboth{$\mQ$-Fano 3-fold and FTS }{$\mQ$-Fano 3-fold and FTS}
{\small{}{\tableofcontents{}}}{\small\par}

2020\textit{ Mathematics subject classification}: 14J45, 14E30, 17A40,17C50


\section{\textbf{Introduction}}

%Should remove quot from Tak3,4 change to Tak6,7

%Gorensteinness $9\times 16$ resolution

\subsection{Classification of $\mQ$-Fano threefolds}

A complex projective variety is called a \textsl{$\mathbb{Q}$-Fano
variety} if it is a normal variety with only terminal singularities,
and its anti-canonical divisor is ample. A $\mathbb{Q}$-Fano variety
is called \textit{prime} if its anti-canonical divisor generates the
group of numerical equivalence classes of $\mQ$-Cartier divisors.
This paper concerns with the classification of prime $\mQ$-Fano 3-folds
and is a companion paper to \cite{Tak4,Tak5,Tak6,Tak7}, where we
construct certain affine varieties and show that they produce, as
weighted complete intersections of their appropriate weighted projectivizations,
several examples of prime $\mathbb{Q}$-Fano 3-folds. We call such
affine varieties or its appropriate weighted projectivizations \textit{key
varieties for $\mQ$-Fano $3$-folds.} The prime $\mQ$-Fano 3-folds
constructed in ibid.~are of anti-canonical codimension 4, where,
by the anti-canonical codimension of a $\mathbb{Q}$-Fano 3-fold $X$,
we mean the codimension of $X$ in the weighted projective space of
the minimal dimension determined by the anti-canonical graded ring
of $X$. 

The affine varieties constructed in ibid.~are related with $\mP^{2}\times\mP^{2}$-fibrations.
In this paper, we construct several affine varieties related with
$\mP^{1}\times\mP^{1}\times\mP^{1}$-fibrations, and produce, as appropriate
weighted complete intersections of their weighted projectivizations,
certain examples of prime $\mathbb{Q}$-Fano 3-folds of anti-canonical
codimension 4. Ahead of our study including ibid., Coughlan and Ducat
did similar work in \cite{CD} defining the $C_{2}$- and $G_{2}^{(4)}$-cluster
varieties and using them as key varieties, where the former is related
with a $\mP^{2}\times\mP^{2}$-fibration, and the latter is related
with a $\mP^{1}\times\mP^{1}\times\mP^{1}$-fibration.

Actually, further ahead of \cite{CD} and ibid., Papadakis \cite{PHD,P2}
constructed via the theory of unprojection more general affine varieties
related with $\mP^{2}\times\mP^{2}$- or $\mP^{1}\times\mP^{1}\times\mP^{1}$-fibrations
with nicknames Tom and Jerry respectively. They seem, however, experimentally
too large as key varieties for prime $\mQ$-Fano $3$-folds; in \cite{CD}
and ibid., to produce prime $\mQ$-Fano 3-fold, we extract appropriate
subvarieties from Papadakis' affine varieties. It is non-trivial what
kind of subvarieties are chosen. The cluster varieties were discovered
in \cite{CD} via mirror symmetry of log Calabi-Yau surfaces. One
purpose of this paper is to reveal that the theory of \textit{Freudenthal
triple system} is a natural framework to define our key varieties
and to describe their equations.

\subsection{Freudenthal triple system}

The first example of a Freudenthal triple system (FTS, in short) is
given by Freudenthal in a series of works \cite{Fr}; the example,
which we denote by $V_{F},$ is a real 56-dimensional representation
of the exceptional group of type $E_{7}$ and is constructed from
a real 27-dimensional exceptional Jordan algebra of type $E_{6}$.
This example $V_{F}$ of an FTS is endowed with a tri-linear product
and is associated with a symplectic form and a quartic form. The exceptional
group of type $E_{7}$ is recovered from these two forms as the group
of linear transformations of $V_{F}$ leaving these two forms invariant. 

In 1960s and 70s, axiomatic definitions of FTS were given by several
researches, which start from a vector space with a symplectic form
and a tri-linear product and then a quartic form defined from them.
Among such researches, we basically follow Ferrar's one \cite{Fe}
in this paper; the definition and the way of investigation of FTS
given by him turn out to be quite suitable to describe our key varieties
conceptually.

Ferrar's idea to investigate the structure of FTS is based upon the
concept of \textit{strictly regular element in FTS }(see Subsection
\ref{subsec:Basics-of-Freudenthal} for the precise definition). He
derives the Peirce decomposition of an FTS from a pair of supplementary
strictly regular elements and then give a good coordinatization of
the FTS. For us, it is important to state the conditions of strictly
regular elements in a way suitable for our purpose (Theorem \ref{thm:stregoneeq}).
Then, specializing to an 8-dimensional FTS, we derive another coordinatization
of the FTS from the Ferrar's Peirce decomposition (Theorem \ref{thm:EqF22}). 

\subsection{Key varieties for prime $\mQ$-Fano threefolds}

For a general 8-dimensional FTS over the complex number field $\mC$,
the locus consisting of strictly regular elements is isomorphic to
the affine cone over $\mP^{1}\times\mP^{1}\times\mP^{1}$ (Corollary
\ref{cor:P1P1P1}). From our coordinatization of an 8-dimensional
FTS mentioned above, we may define an affine scheme $\mathfrak{F}_{\mA}^{22}$
related with a $\mP^{1}\times\mP^{1}\times\mP^{1}$-fibration (Definition
\ref{def:In-the-affineF22}). It is a subscheme of Papadakis' Jerry
as mentioned above but is still too large to produce prime $\mQ$-Fano
3-folds. We extract several subvarieties of $\mathfrak{F}_{\mA}^{22}$
and obtain several results as for key varieties of prime $\mQ$-Fano
3-folds. Among other things, we state the following as the main result
of this paper, which summarizes Propositions \ref{prop:The-singular-locus},
\ref{prop:Jerry14}, \ref{prop:UTerm}, and Theorem \ref{thm:20544Ex}:
\begin{thm}
\label{thm:A-variety-14Sum} The following assertions hold:

\begin{enumerate}

\item There exists a $14$-dimensional factorial affine variety $\mathfrak{U}_{\mA}^{14}$
of codimension $4$ in an affine $18$-space with only Gorenstein
terminal singularities. 

\item Examples of $\mQ$-Fano $3$-folds of No.20544 in \cite{GRDB}
are produced as weighted complete intersections of the weighted projectivization
of $\mathfrak{U}_{\mA}^{14}$ in the weighted projective space $\mP(1^{15},2^{2},3)$. 

\end{enumerate}
\end{thm}

We remark that examples of $\mQ$-Fano $3$-folds of No.20544 are
constructed from the $C_{2}$-cluster variety but are not constructed
from the $G_{2}^{(4)}$-cluster variety (\cite{bigtables}). 

\subsection{Structure of the paper\label{subsec:Structure-of-the}}

The paper consists of two parts.

In Part 1, we revisit the theory of FTS in detail basically following
\cite{Fe} with emphasis on strictly regular elements (in some places,
we also refer to \cite{B} and \cite{Kr}). The main result of Part
1 is Theorem \ref{thm:stregoneeq}, which states the conditions of
strictly regular elements and leads to the equations of key varieties
given in Part 2. 

In \cite{Fe}, it is shown that an FTS with a nondegenerate skew form
has a direct sum decomposition with two copies of a Jordan algebra
of a cubic form as direct summands. For our purpose, however, we prefer
not to identify Jordan algebras in FTS fully though several concepts
for FTS are helpfully understood by Jordan algebraic considerations.
For this reason, we decide to revisit the theory of FTS in detail
in this paper.

In Part 2, we first construct an affine scheme $\mathfrak{F}_{\mathbb{A}}^{22}$
from a natural parameterization of an 8-dimensional FTS. In Sections
\ref{sec:KeyEquation14}--\ref{sec:Z12}, we extract three subvarieties
$\mathfrak{U}_{\mathbb{A}}^{14}$ ,$\mathfrak{S}_{\mathbb{A}}^{8}$,
$\mathfrak{Z}_{\mathbb{A}}^{12}$ respectively of $\mathfrak{F}_{\mathbb{A}}^{22}$
and show that they produce, as weighted complete intersections of
their weighted projectivizations, certain prime $\mQ$-Fano 3-folds.
Several properties of $\mathfrak{U}_{\mathbb{A}}^{14}$ ,$\mathfrak{S}_{\mathbb{A}}^{8}$,
$\mathfrak{Z}_{\mathbb{A}}^{12}$ are also obtained and they are important
to construct prime $\mQ$-Fano 3-folds. We describe $\mathfrak{U}_{\mathbb{A}}^{14}$
in detail in Section \ref{sec:KeyEquation14} (we refer to Theorem
\ref{thm:A-variety-14Sum} for a summary). As is clarified in this
paper, the affine variety $\mathfrak{U}_{\mathbb{A}}^{14}$ is also
the cornerstone for studying other key varieties. Moreover, it will
produce other prime $\mQ$-Fano 3-folds in our future work. By these
reasons, we devote many pages to the investigation of $\mathfrak{U}_{\mathbb{A}}^{14}$.
As for the variety $\mathfrak{S}_{\mA}^{8}$, we describe in detail
its equations, an ${\rm SL}_{2}\rtimes{\rm SL}_{2}$-action on it,
and its relation with a $\mP^{1}\times\mP^{1}\times\mP^{1}$-fibration.
We revisit $\mathfrak{S}_{\mA}^{8}$ in Section \ref{sec:The--Cluster-variety}.
As for the variety $\mathfrak{Z}_{\mathbb{A}}^{12}$, we show in Section
\ref{sec:Z12} that a certain weighted projectivization of it is the
$11$-dimensional $\mQ$-Fano variety constructed in the paper \cite{Tak3}. 

In Section \ref{sec:The--Cluster-variety}, we revisit the $G_{2}^{(4)}$-cluster
variety defined in \cite{CD}, which we denote by $\mathfrak{Cl}_{\mA}^{10}$.
By the big table \cite{bigtables} of prime $\mQ$-Fano 3-folds obtained
from $\mathfrak{Cl}_{\mA}^{10}$, we observe that such prime $\mQ$-Fano
3-folds are weighted complete intersections of some weighted projectivizations
of $\mathfrak{Cl}_{\mA}^{10}$ itself or its several subvarieties.
In Subsections \ref{subsec:The-case-onlyone}--\ref{subsec:The-case-withtwoparam.},
we show that $\mathfrak{Cl}_{\mA}^{10}$ itself and such subvarieties
are actually isomorphic to subvarieties of $\mathfrak{U}_{\mathbb{A}}^{14}$
or $\mathfrak{Z}_{\mathbb{A}}^{12}$ weighted homogeneously with respect
to some weights on the coordinates. Hence it turns out that all prime
$\mQ$-Fano 3-folds obtained from $\mathfrak{Cl}_{\mA}^{10}$ are
also obtained from $\mathfrak{U}_{\mathbb{A}}^{14}$ or $\mathfrak{Z}_{\mathbb{A}}^{12}$.
Among other things, we also show that the subvariety of $\mathfrak{Cl}_{\mA}^{10}$
studied in Subsection \ref{subsec:The-case-A3A4} is weighted homogeneously
isomorphic to $\mathfrak{S}_{\mA}^{8}$, and the subvarieties $\mathfrak{T}_{\mA}^{8}$
and $\mathfrak{B}_{\mA}^{6}$ of $\mathfrak{Cl}_{\mA}^{10}$ studied
in Subsections \ref{subsec:The-case-AA1A4} and \ref{subsec:The-case-withtwoparam.}
admit a nontrivial ${\rm SL}_{2}$-action. 

Finally, in Section \ref{sec:Affine-variety-Papadakis}, we see another
aspect of the variety $\mathfrak{F}_{\mA}^{22}$; we show that an
open subset of the $23$-dimensional affine variety $\mathfrak{P}_{\mA}^{23}$
constructed by Papadakis using Type ${\rm II}_{1}$ unprojection \cite{P3,P4}
is transformed to an open subset of the cone over $\mathfrak{F}_{\mA}^{22}$. 

\vspace{5pt}

This paper contains several assertions which can be proved by straightforward
computations; we often omit such computations. Some computations are
difficult by hand but are easy within a software package. In our computations,
we use intensively the software systems Mathematica \cite{W} and
\textsc{Singular} \cite{DGPS}. 

\subsection{A future plan}

As is noted in Subsection \ref{subsec:Structure-of-the}, the affine
variety $\mathfrak{U}_{\mathbb{A}}^{14}$ produces examples of prime
$\mQ$-Fano $3$-folds which are also obtained from the $G_{2}^{(4)}$-cluster
variety $\mathfrak{Cl}_{\mA}^{10}$. By Theorem \ref{thm:A-variety-14Sum}
(2), one more example is already added in this paper. Actually, we
can verify that the affine variety $\mathfrak{U}_{\mA}^{14}$ produce
more examples of prime $\mQ$-Fano 3-folds. Since the verification
takes more pages (cf.~\cite{Tak7}), we will publish it elsewhere.

\vspace{5pt}

\noindent\textbf{ Notation.}
\begin{itemize}
\item $w(*)$: the weight on the coordinate $*$ of a weighted projective
space.
\item If we put weights for the entries of a matrix $A$, we denote by $w(A)$
the set of the weights implemented in the matrix form corresponding
to $A$.
\item For a $2\times3$ matrix $M=\left(\begin{array}{ccc}
f_{1} & f_{2} & f_{3}\\
g_{1} & g_{2} & g_{3}
\end{array}\right)$, we define \begin{equation*} \wedge^2 M=\begin{pmatrix}\begin{vmatrix}f_2 & f_3 \\ g_2 & g_3\end{vmatrix}&-\begin{vmatrix}f_1 & f_3 \\ g_1 & g_3\end{vmatrix}& \begin{vmatrix}f_1 & f_2 \\ g_1 & g_2\end{vmatrix} \end{pmatrix}.\end{equation*}
For a $3\times2$ matrix, we have a similar definition. 
\end{itemize}
\vspace{5pt}

\noindent\textbf{ Acknowledgment}. I am grateful to Professor Shigeru
Mukai; inspired by his articles \cite{Mu1,Mu2}, I was led to Jordan
algebras and FTSs to describe key varieties. Moreover, I was helped
to find the key variety $\mathfrak{U}_{\mA}^{14}$ by his inference
for the dimensions of key varieties. I would like to dedicate this
paper to him on his 70th birthday. This work is supported in part
by Grant-in Aid for Scientific Research (C) 16K05090.

\part{Freudenthal triple systems\label{part:Freudenthal-triple-systems}}

\section{\textbf{Strictly regular elements in Freudenthal triple systems\label{sec:Strictly-regular-elements}}}

In this section, we basically follow the treatment of Freudenthal
triple system by Ferrar \cite{Fe} with modifications in several coefficients
of equalities after Brown \cite{B}. 

Throughout this section, we assume that $\mathsf{k}$ is a field of
characteristic $\not=2,3$.

\subsection{Basics of Freudenthal triple system\label{subsec:Basics-of-Freudenthal}}
\begin{defn}[{\cite[Sec.1]{Fe}}]
\label{def:Triple} A \textit{Freudenthal triple system (FTS for
short)} is a $\mathsf{k}$-vector space $V$ with a tri-linear product
\[
V\times V\times V\ni(p_{1},p_{2},p_{3})\mapsto p_{1}\bullet p_{2}\bullet p_{3}\in V
\]
and a skew bi-linear form 
\[
V\times V\ni(p_{1},p_{2})\mapsto\omega(p_{1},p_{2})\in\mathsf{k}
\]
such that 

\begin{enumerate}[({A}1)]

\item the tri-linear product is symmetric in all arguments,

\item
\[
\widetilde{F}(p_{1},p_{2},p_{3},p_{4}):=\omega(p_{1}\bullet p_{2}\bullet p_{3},p_{4})
\]

is a nonzero symmetric $4$-linear form, and

\item the equality
\[
6(p\bullet p\bullet p)\bullet p\bullet q=\omega(q,p)(p\bullet p\bullet p)+\omega\left(q,p\bullet p\bullet p\right)p
\]
holds for any $p,q\in V$.

\end{enumerate}
\end{defn}

\begin{rem}
We note the coefficient $6$ in the l.h.s.~of the equality in (A3).
This formulation is according to Brown \cite{B} and this influences
several coefficients in the equalities below although we will not
mention one by one. 
\end{rem}

Linearizing the equality in (A3) completely, we obtain the following:
\begin{align}
\label{eq:penta}\\
 & 6\left((p_{1}\bullet p_{2}\bullet p_{3})\bullet p_{4}+(p_{1}\bullet p_{2}\bullet p_{4})\bullet p_{3}+(p_{1}\bullet p_{3}\bullet p_{4})\bullet p_{2}+(p_{2}\bullet p_{3}\bullet p_{4})\bullet p_{1}\right)\bullet q=\nonumber \\
 & \omega(q,p_{4})(p_{1}\bullet p_{2}\bullet p_{3})+\omega(q,p_{3})(p_{1}\bullet p_{2}\bullet p_{4})+\omega(q,p_{2})(p_{1}\bullet p_{3}\bullet p_{4})+\omega(q,p_{1})(p_{2}\bullet p_{3}\bullet p_{4})\nonumber \\
+ & \omega\left(q,p_{1}\bullet p_{2}\bullet p_{3}\right)p_{4}+\omega\left(q,p_{1}\bullet p_{2}\bullet p_{4}\right)p_{3}+\omega\left(q,p_{1}\bullet p_{3}\bullet p_{4}\right)p_{2}+\omega\left(q,p_{2}\bullet p_{3}\bullet p_{4}\right)p_{1}.\nonumber 
\end{align}

In this paper, we call this \textit{the pentagram product formula.}

For $p\in V,$ we denote by 
\[
L_{p,p}\colon V\to V
\]
 the linear map defined by
\[
V\ni q\mapsto p\bullet p\bullet q\in V.
\]

The following concept plays a central role in this paper. 
\begin{defn}[{\cite[Sec.3]{Fe}}]
 An element $p\in V$ is called \textit{a strictly regular element}
if it holds that 
\[
{\rm Image}L_{p,p}\subset\mathsf{k}p.
\]
\end{defn}

\begin{prop}[{\cite[p.317, (5) and Lem.3.1]{Fe}}]
\label{prop:streg} An element $p\in V$ is strictly regular if and
only if it holds that 
\[
3L_{p,p}(q)+\omega(p,q)p=0
\]
 for any $q\in V$.
\end{prop}

\subsection{Jordan algebraic description of FTS\label{subsec:Jordan-algebraic-description}}

In \cite{Fe}, it is shown that an FTS with nondegenerate $\omega$
has a direct sum decomposition with two copies of a Jordan algebra
of cubic form as direct summands. In this subsection, we proceed slightly
in a different way without identifying Jordan algebra fully, which
is suitable for our purpose. 
\begin{prop}[{\cite[Sec.4 (6)]{Fe}}]
 Let $e_{s},e_{t}\in V$ be supplementary, strictly regular elements
(namely, $e_{s},e_{t}$ are strictly regular such that $\omega(e_{s},e_{t})=1$).
The following equality holds:
\[
L_{e_{s},e_{t}}^{2}p=1/12\,\omega(p,e_{t})e_{s}-1/12\,\omega(p,e_{t})e_{t}+1/36\,p.
\]
\end{prop}

Hereafter in this section, we assume that \textit{
\begin{equation}
\text{the skew bilinear form \ensuremath{\omega} is nondegenerate.}\label{eq:nondeg}
\end{equation}
}
\begin{cor}[{\cite[Sec.4]{Fe}}]
 Let $V_{\alpha}$ be the $L_{e_{s},e_{t}}$-eigenspace for the eigenvalue
$\alpha$. The vector space $V$ has the following decomposition into
the $L_{e_{s},e_{t}}$-eigenspaces:
\[
V=V_{-1/3}\oplus V_{1/3}\oplus V_{1/6}\oplus V_{-1/6},
\]
where $V_{-1/3}=\mathsf{k}e_{s}$ and $V_{1/3}=\mathsf{k}e_{t}$.
This decomposition is called \textup{the} \textup{Peirce decomposition}. 
\end{cor}

Hereafter we use the following notation:
\[
V_{s}:=V_{-1/3},\,V_{t}:=V_{1/3},\,V_{x}:=V_{1/6},\,V_{y}:=V_{-1/6}.
\]

\vspace{5pt}

Now we begin Jordan algebraic treatments of FTS.\textit{ }We refer
to \cite[Sec.4.2]{Mc2} for the subjects in the theory of Jordan algebra
corresponding to those in the sequel. 

We set 

\begin{equation}
\text{\ensuremath{\beta(x,y):=\omega(x,y)} for \ensuremath{x\in V_{x}} and \ensuremath{y\in V_{y}}.}\label{eq:trace}
\end{equation}
 The bi-linear form $\beta$ corresponds to that for a Jordan algebra
of a cubic form \cite[p.189]{Mc2}. As noted in \cite[p.318, the 3rd line from the bottom]{Fe},
\begin{equation}
\text{\ensuremath{\beta(x,y)} is non-degenerate.}\label{eq:betanondeg}
\end{equation}
Hence $\dim V_{x}=\dim V_{y}$, which we will denote by $n$;
\[
n:=\dim V_{x}=\dim V_{y}.
\]
For $x\in V_{x}$ and $y\in V_{y}$, we set 
\begin{equation}
N_{x}(x):=1/2\,\omega(e_{s},x\bullet x\bullet x),\,N_{y}(y):=1/2\,\omega(e_{t},y\bullet y\bullet y),\label{eq:Norm}
\end{equation}

and 
\begin{equation}
x^{\sharp}:=3/2\,x\bullet x\bullet e_{s},\,y^{\sharp}:=-3/2\,y\bullet y\bullet e_{t}.\label{eq:Ad}
\end{equation}
By \cite[Lem.4.1]{Fe}, we have $x^{\sharp}\in V_{y}$ and $y^{\sharp}\in V_{x}$.
The cubic forms $N_{x}$ and $N_{y}$ correspond to the cubic form
for a Jordan algebra of a cubic form.

For $x_{1},x_{2}\in V_{x}$ and $y_{1},y_{2}\in V_{y}$, we also set
\begin{equation}
x_{1}\sharp x_{2}:=(x_{1}+x_{2})^{\sharp}-x_{1}^{\sharp}-x_{2}^{\sharp},\,y_{1}\sharp y_{2}:=(y_{1}+y_{2})^{\sharp}-y_{1}^{\sharp}-y_{2}^{\sharp}.\label{eq:bisha}
\end{equation}
Then, by (\ref{eq:Ad}), we have 
\begin{equation}
x_{1}\sharp x_{2}=3x_{1}\bullet x_{2}\bullet e_{s},\,y_{1}\sharp y_{2}=-3y_{1}\bullet y_{2}\bullet e_{t}.\label{eq:sharp}
\end{equation}

\begin{prop}
\label{prop:NxNy} For $x\in V_{x}$ and $y\in V_{y},$ it holds that
\begin{equation}
N_{x}(x)=1/3\,\beta\left(x,x^{\sharp}\right),N_{y}(y)=1/3\,\beta\left(y^{\sharp},y\right).\label{eq:normxy}
\end{equation}
\end{prop}

\begin{proof}
By (\ref{eq:Norm}), (A2) in Definition \ref{def:Triple} and (\ref{eq:Ad}),
we have the first equality as follows: 
\[
N_{x}(x)=1/2\omega(e_{s},x\bullet x\bullet x)=1/2\omega(x,x\bullet x\bullet e_{s})=1/2\beta(x,2/3x^{\sharp})=1/3\beta(x,x^{\sharp}).
\]

The second one follows similarly.
\end{proof}
\begin{cor}
\label{cor:partialN}For $x,x'\in V_{x}$ and $y,y'\in V_{y},$ it
holds that 
\begin{equation}
\partial_{x'}N_{x}(x)=\beta\left(x',x^{\sharp}\right),\partial_{y'}N_{y}(y)=\beta\left(y^{\sharp},y'\right),\label{eq:normxypartial}
\end{equation}
where $\partial_{x'}N_{x}(x)$ is the directional derivative of $N_{x}$
in the direction $x'$, evaluated at $x$, and $\partial_{y'}N_{y}(y)$
is similarly defined.
\end{cor}

\begin{proof}
By (\ref{eq:Norm}), we have $\partial_{x'}N_{x}(x)=1/2\,\omega(e_{s},3x\bullet x\bullet x'),$
and the r.h.s. is equal to $3/2\,\omega(x',x\bullet x\bullet e_{s})=3/2\omega(x',2/3x^{\sharp})=\beta(x',x^{\sharp})$
by (A2) in Definition \ref{def:Triple} and (\ref{eq:Ad}). Therefore
the first equality follows. The second one follows similarly.
\end{proof}
Corollary \ref{cor:partialN} corresponds to Trace-Sharp formula in
\cite[p.189]{Mc2}.

\vspace{3pt}

By \cite{Fe}, we have 
\begin{lem}
The following equalities hold:
\begin{equation}
(x\bullet x\bullet e_{s})\bullet(x\bullet x\bullet e_{s})\bullet e_{t}=4/27\omega(x\bullet x\bullet x,e_{s})x\label{eq:xxsxxst}
\end{equation}
\[
x\bullet x\bullet y=-1/3\,\omega(x,y)x-3(x\bullet x\bullet e_{s})\bullet e_{t}\bullet y
\]
\end{lem}

\begin{prop}
\label{prop:=000023=000023} For $x\in V_{x}$ and $y\in V_{y},$
it holds that 
\[
(x^{\sharp})^{\sharp}=N_{x}(x)x,~(y^{\sharp})^{\sharp}=N_{y}(y)y.
\]
\end{prop}

\begin{proof}
By (\ref{eq:Ad}), (\ref{eq:xxsxxst}) and (\ref{eq:Norm}), we have
the first equality as follows:
\begin{align*}
(x^{\sharp})^{\sharp} & =(3/2\,x\bullet x\bullet e_{s})^{\sharp}=(3/2)^{2}\left(-3/2\left(x\bullet x\bullet e_{s}\right)\bullet\left(x\bullet x\bullet e_{s}\right)\bullet e_{t}\right)\\
= & -(3/2)^{3}\left(4/27\,\omega(x\bullet x\bullet x,e_{s})x\right)=-1/2\,\omega(x\bullet x\bullet x,e_{s})x=N_{x}x.
\end{align*}

The second one follows similarly.
\end{proof}
Proposition \ref{prop:=000023=000023} corresponds to Adjoint Identity
in \cite[p.189]{Mc2}.

\vspace{3pt}

The following two auxiliary results are needed to show Theorem \ref{thm:EqF22}.
\begin{cor}
\label{cor:no multfacN} Assume that $n\geq3$ and $N_{x}(x)$ is
not identically zero. Then $N_{x}(x)$ has no multiple factors. The
similar statement for $N_{y}(y)$ also holds.
\end{cor}

\begin{proof}
Assume that $N_{x}(x)$ has a multiple factor $l(x)$, which must
be a linear form since $N_{x}(x)$ is a cubic form. We choose coordinates
$x_{1},\dots,x_{n}$ of $V_{x}$ such that $l(x)=x_{1}$ and $N_{x}(x)=x_{1}^{2}x_{2}$
(recall that $n\geq2).$ Then, for $a={\empty}^{t}\!(a_{1},\dots,a_{n}),$
$b={\empty}^{t}\!(b_{1},\dots,b_{n})\in V_{x}$ (we consider $a,b$
as column vectors), we have $\partial_{a}N_{x}(b)=2a_{1}b_{1}b_{2}+a_{2}b_{1}^{2}$,
which is equal to $\beta(a,b^{\sharp})$ by Corollary \ref{cor:partialN}.
Note that the subset $\{b^{\sharp}\mid b\in V_{x}\}\subset V_{y}$
is dense in $V_{y}$ by Proposition \ref{prop:=000023=000023} and
the assumption that $N_{x}(x)$ is not identically zero. Therefore,
by the assumption that $n\geq3$, the equality $\beta(a,b^{\sharp})=2a_{1}b_{1}b_{2}+a_{2}b_{1}^{2}$
implies that $\beta$ is degenerate , a contradiction to (\ref{eq:betanondeg}).

We can show the claim for $N_{y}(y)$ similarly.
\end{proof}
\begin{cor}
\label{cor:nocommonfac=000023} Assume that $n\geq3$ and $N_{x}(x)$
is not identically zero. The coordinates of $x^{\sharp}$ have no
common factors and the similar statement holds for $y^{\sharp}$.
\end{cor}

\begin{proof}
If the coordinates of $x^{\sharp}$ have a common factor $F(x)$,
we see that $F(x)^{2}$ divides $N_{x}(x)$ by Proposition \ref{prop:=000023=000023}
since the operation $\sharp$ is quadratic. This contradicts Corollary
\ref{cor:no multfacN}. We can show the claim for $y^{\sharp}$ similarly.
\end{proof}
The following formula is useful for calculations below.
\begin{lem}
\label{lem:xxsha}For $x,\widetilde{x}\in V_{x}$ and $y\in V_{y}$,
it holds that 
\begin{equation}
x\bullet x^{\sharp}\bullet\widetilde{x}=1/6\,(\beta(\widetilde{x},x^{\sharp})x-N_{x}(x)\widetilde{x}),\label{eq:xxshax'}
\end{equation}
\begin{equation}
x\bullet x^{\sharp}\bullet y=1/6\,(-\beta(x,y)x^{\sharp}+N_{x}(x)y),\label{eq:xxshay}
\end{equation}
 and 
\begin{equation}
y\bullet y^{\sharp}\bullet x=1/6\,(\beta(x,y)y^{\sharp}-N_{y}(y)x).\label{eq:yyshax}
\end{equation}
\end{lem}

\begin{proof}
We only give a proof of the formula (\ref{eq:xxshax'}) since we can
show the remaining two similarly. By (\ref{eq:Ad}), we have $x\bullet x^{\sharp}\bullet\widetilde{x}=3/2(x\bullet x\bullet e_{s})\bullet x\bullet\widetilde{x}$.
Applying the pentagram product formula (\ref{eq:penta}) with (\ref{eq:Ad})
and (\ref{eq:normxy}), we obtain (\ref{eq:xxshax'}).
\end{proof}
Hereafter we also follow the paper \cite{Kr} with the above treatment
in this subsection.

\vspace{3pt}

\noindent \textbf{Cube Formula (\cite[p.946,(48)]{Kr}). }For $p=(s,t,x,y)\in V=V_{s}\oplus V_{t}\oplus V_{x}\oplus V_{y},$
it holds that\textbf{ 
\begin{align}
p\bullet p\bullet p & =\label{eq:ppp}\\
 & (-s^{2}t+s\beta(x,y)-2N_{y})e_{s}+(st^{2}-t\beta(x,y)+2N_{x})e_{t}\nonumber \\
 & +\left(st-\beta(x,y)\right)x+2x^{\sharp}\sharp y-2ty^{\sharp}\nonumber \\
 & -\left(st-\beta(x,y)\right)y-2x\sharp y^{\sharp}+2sx^{\sharp}.\nonumber 
\end{align}
}
\begin{proof}
In \cite{Kr}, the proof is given for an FTS which contains two copies
of a Jordan algebra of a cubic form. In any way, we may verify the
assertion by a straightforward calculation using the facts obtained
so far. 
\end{proof}

\subsection{Condition on strict regularity}

Via partial linearization of (\ref{eq:ppp}) (with a minor correction),
we arrive at the following condition on strict regularity by Proposition
\ref{prop:streg}: 
\begin{prop}[{\textbf{Condition on strict regularity (\cite[p.947,(53)]{Kr})}}]
\label{prop:cond streg} 
\end{prop}

For $p=(s,t,x,y)$, $q=(\widetilde{s},\widetilde{t,}\widetilde{x},\widetilde{y})\in V$,
it holds that

\begin{align*}
 & 3p\bullet p\bullet q+\omega(p,q)p=\\
 & \left(-\widetilde{s}\left(3st-\beta(x,y)\right)+2\beta(sx-y^{\sharp},\widetilde{y})\right)e_{s}\\
 & +\left(\left(3st-\beta(x,y)\right)\widetilde{t}-2\beta(\widetilde{x},ty-x^{\sharp})\right)e_{t}\\
 & +\left(st-1/3\,\beta(x,y)\right)\widetilde{x}+2\widetilde{t}\left(sx-y^{\sharp}\right)-2\widetilde{y}\sharp(ty-x^{\sharp})+2\Delta_{x}(x,y;\widetilde{x})\\
 & -\left(st-1/3\,\beta(x,y)\right)\widetilde{y}-2\widetilde{s}\left(ty-x^{\sharp}\right)+2\widetilde{x}\sharp(sx-y^{\sharp})-2\Delta_{y}(x,y;\widetilde{y}),
\end{align*}
where we set 
\begin{align*}
 & \Delta_{x}(x,y;\widetilde{x}):=-1/3\,\beta(x,y)\widetilde{x}-\beta(\widetilde{x},y)x+(x\sharp\widetilde{x})\sharp y,\\
 & \Delta_{y}(x,y;\widetilde{y}):=-1/3\,\beta(x,y)\widetilde{y}-\beta(x,\widetilde{y})y+(y\sharp\widetilde{y})\sharp x.
\end{align*}
Computing $(x\sharp\widetilde{x})\sharp y$ by the pentagram product
formula (\ref{eq:penta}), we see that 

\begin{equation}
\Delta_{x}(x,y;\widetilde{x})=1/6\,\beta(x,y)\widetilde{x}-1/2\,\beta(\widetilde{x},y)x+3x\bullet\widetilde{x}\bullet y.\label{eq:Dxtrans}
\end{equation}
Similarly, we have 
\begin{equation}
\Delta_{y}(x,y;\widetilde{y})=1/6\,\beta(x,y)\widetilde{y}-1/2\,\beta(x,\widetilde{y})y-3y\bullet\widetilde{y}\bullet x.\label{eq:Dytrans}
\end{equation}

By Propositions \ref{prop:streg} and \ref{prop:cond streg}, we have
the following:
\begin{cor}
\label{cor:streg}An element $p=(s,t,x,y)\in V$ is strictly regular
if and only if it holds that 
\[
sx=y^{\sharp},\,ty=x^{\sharp},\,st=1/3\,\beta(x,y),
\]
and
\begin{equation}
\Delta_{x}(x,y;\widetilde{x})=\Delta_{y}(x,y;\widetilde{y})=0\,\,\text{for any \ensuremath{\widetilde{x}\in V_{x}},\,\ensuremath{\widetilde{y}\in V_{y}}}.\label{eq:forany}
\end{equation}
\end{cor}

\vspace{3pt}

Now we will see that only the conditions on $\Delta_{x}$ in (\ref{eq:forany})
is necessary. This is inspired by \cite[Lem.23]{Kr}. We give here
a proof in a different flavor from the argument there. 
\begin{lem}
\label{lem:DxDy} For $x,\widetilde{x}\in V_{x}$ and $y,\widetilde{y}\in V_{y}$,
the following two equalities hold:

\begin{equation}
\beta\left(\Delta_{x}(x,y;\widetilde{x}),\widetilde{y}\right)=\beta\left(\widetilde{x},\Delta_{y}(x,y;\widetilde{y})\right).\label{eq:DxDy1}
\end{equation}

\begin{equation}
\beta\left(\Delta_{x}(x,y;\widetilde{x}),\widetilde{y}\right)=\beta\left(x,\Delta_{y}(\widetilde{x},\widetilde{y};y)\right).\label{eq:DxDy2}
\end{equation}
\end{lem}

\begin{proof}
The equation (\ref{eq:DxDy1}) follows from the following chain of
the equalities:
\begin{align*}
 & \beta\left(\Delta_{x}(x,y;\widetilde{x}),\widetilde{y}\right)\\
 & =\beta\left(1/6\,\beta(x,y)\widetilde{x}-1/2\,\beta(\widetilde{x},y)x+3x\bullet\widetilde{x}\bullet y,\widetilde{y}\right)(\text{by}\,\,(\ref{eq:Dxtrans}))\\
 & =1/6\,\beta(x,y)\beta(\widetilde{x},\widetilde{y})-1/2\,\beta(\widetilde{x},y)\beta(x,\widetilde{y})+\beta(3x\bullet\widetilde{x}\bullet y,\widetilde{y)}\\
 & =1/6\,\beta(x,y)\beta(\widetilde{x},\widetilde{y})-1/2\,\beta(\widetilde{x},y)\beta(x,\widetilde{y})+\omega(3x\bullet\widetilde{y}\bullet y,\widetilde{x})\,(\text{by (A2) in Def.\,\ref{def:Triple}})\\
 & =1/6\,\beta(x,y)\beta(\widetilde{x},\widetilde{y})-1/2\,\beta(\widetilde{x},y)\beta(x,\widetilde{y})-\beta(\widetilde{x},3x\bullet\widetilde{y}\bullet y)\\
 & =\beta\left(\widetilde{x},1/6\,\beta(x,y)\widetilde{y}-1/2\,\beta(x,\widetilde{y})y-3x\bullet\widetilde{y}\bullet y\right)\\
 & =\beta\left(\widetilde{x},\Delta_{y}(x,y;\widetilde{y})\right)(\text{by}\,\,(\ref{eq:Dytrans})).
\end{align*}
To obtain the equation (\ref{eq:DxDy2}), we have only to change the
last 4 lines of the above chain of the equalities as follows: 
\begin{align*}
 & =1/6\,\beta(x,y)\beta(\widetilde{x},\widetilde{y})-1/2\,\beta(\widetilde{x},y)\beta(x,\widetilde{y})+\omega(3\widetilde{x}\bullet\widetilde{y}\bullet y,x)\,(\text{by (A2) in Def.\,\ref{def:Triple}})\\
 & =1/6\,\beta(x,y)\beta(\widetilde{x},\widetilde{y})-1/2\,\beta(\widetilde{x},y)\beta(x,\widetilde{y})-\beta(x,3\widetilde{x}\bullet\widetilde{y}\bullet y)\\
 & =\beta\left(x,1/6\,\beta(\widetilde{x},\widetilde{y})y-1/2\,\beta(\widetilde{x},y)\widetilde{y}-3\widetilde{x}\bullet\widetilde{y}\bullet y\right)\\
 & =\beta\left(x,\Delta_{y}(\widetilde{x},\widetilde{y};y)\right)(\text{by}\,\,(\ref{eq:Dytrans})).
\end{align*}
\end{proof}
\begin{cor}
\label{cor:The-following-conditionsDelxDely}The following conditions
on $x\in V_{x}$ and $y\in V_{y}$ are equivalent:

\begin{enumerate}

\item $\Delta_{x}(x,y;\widetilde{x})=0\,\,\text{for any \ensuremath{\widetilde{x}\in V_{x}}.}$

\item $\Delta_{y}(x,y;\widetilde{y})=0\,\,\text{for any \ensuremath{\widetilde{y}\in V_{y}}}$.

\end{enumerate}
\end{cor}

\begin{proof}
Assume that $\Delta_{x}(x,y;\widetilde{x})=0\,\,\text{for any \ensuremath{\widetilde{x}\in V_{x}}.}$
Then, by the equality (\ref{eq:DxDy1}) as in Lemma \ref{lem:DxDy},
it holds that $\beta\left(\widetilde{x},\Delta_{y}(x,y;\widetilde{y})\right)=0$
for any $\widetilde{x}\in V_{x}$ and $\widetilde{y}\in V_{y}$. Since
$\beta$ is nondegenerate, $\Delta_{y}(x,y;\widetilde{y})=0\,\,\text{for any \ensuremath{\widetilde{y}\in V_{y}}}$.
Thus (1) implies (2). Similarly we see that (2) implies (1).
\end{proof}
Now we arrive at the main result of Part \ref{part:Freudenthal-triple-systems}
by Corollaries \ref{cor:streg} and \ref{cor:The-following-conditionsDelxDely}.
\begin{thm}
\label{thm:stregoneeq} An element $p=(s,t,x,y)$ is strictly regular
if and only if it holds that 
\begin{align}
sx & =y^{\sharp},\,ty=x^{\sharp},\label{eq:sx ty}
\end{align}
\begin{equation}
st=1/3\,\beta(x,y),\label{eq:st}
\end{equation}
and
\begin{equation}
\Delta_{x}(x,y,\widetilde{x})=0\,\,\text{for any }\text{\ensuremath{\widetilde{x}}\ensuremath{\in V_{x}}}.\label{eq:forany2}
\end{equation}
\end{thm}

\begin{defn}
\label{def:We-denote-bystregR}We denote by $\mathfrak{R}$ the affine
scheme in $V$ defined by the equations (\ref{eq:sx ty}), (\ref{eq:st}),
and (\ref{eq:forany2}) as in Theorem \ref{thm:stregoneeq}.
\end{defn}

\begin{rem}
Here we mention some background of Theorem \ref{thm:stregoneeq}.
Our derivation of the defining equations of $\mathfrak{R}$ as in
Theorem \ref{thm:stregoneeq} is inspired by \cite[Lem.23]{Kr} and
\cite[p.253--254]{YamA}, and its origin is traced back to the series
of the fundamental papers \cite{Fr} (see also \cite[Prop.6.2]{Cl}
and \cite[(2.1)]{Fa}). We refer to \cite[Prop.6.2]{Cl} and \cite[p.515, Sec.0]{KaYas}
for aspects of this equation of $\mathfrak{R}$ in a graded Lie algebra
of contact type. In \cite{KaYas}, the projectivization of $\mathfrak{R}$
is called Freudenthal variety. In \cite[Sec.5]{Mu1} and \cite{Mu2},
the projectivization of $\mathfrak{R}$ is studied in relation with
smooth Fano threefolds and is called Legendre projective variety. 
\end{rem}

\subsection{More on $\Delta_{x}$ and $\Delta_{y}$}

In this subsection, we examine the condition (\ref{eq:forany2}) more
in detail. 

By (\ref{eq:normxy}), Lemma \ref{lem:xxsha}, (\ref{eq:Dxtrans}),
and (\ref{eq:Dytrans}), we immediately obtain the following:
\begin{lem}
\label{lem:xshaId}For any $x,\widetilde{x}\in V_{x}$ and any $y,\widetilde{y}\in V_{y}$,
it holds that

\[
\Delta_{x}(x,x^{\sharp},\widetilde{x})=0,\Delta_{x}(y^{\sharp},y,\widetilde{x})=0,
\]
and 

\[
\Delta_{y}(y^{\sharp},y,\widetilde{y})=0,\Delta_{y}(x,x^{\sharp},\widetilde{y})=0.
\]
\end{lem}

\begin{lem}
\label{lem:n-1dim}The following two conditions for $\widetilde{x}\in V_{x}$
with $N_{x}(\widetilde{x})\not=0$ and $\widetilde{y}\in V_{y}$ are
equivalent:

\begin{enumerate}

\item It holds that 
\begin{equation}
\beta\left(\Delta_{x}(x,y;\widetilde{x}),\widetilde{y}\right)=0\label{eq:csharp}
\end{equation}
 for any $x\in V_{x}$ and $y\in V_{y}$.

\item $\widetilde{y}\in k\cdot\widetilde{x}^{\sharp}$.

\end{enumerate}
\end{lem}

\begin{proof}
By (\ref{eq:DxDy2}) as in Lemma \ref{lem:DxDy}, we have $\beta\left(\Delta_{x}(x,y;\widetilde{x}),\widetilde{y}\right)=\beta\left(x,\Delta_{y}(\widetilde{x},y;\widetilde{y})\right)$.
Therefore, (1) is equivalent to that $\beta\left(x,\Delta_{y}(\widetilde{x},y;\widetilde{y})\right)=0$
for any $x\in V_{x}$ and $y\in V_{y}$. By nondegeneracy of $\beta$,
this is also equivalent to 
\begin{equation}
\Delta_{y}(\widetilde{x},\widetilde{y};y)=0\label{eq:Dy}
\end{equation}
for any $y\in V_{y}$.

\vspace{3pt}

We show (1) implies (2). Assume that (1) holds. We set $y=\widetilde{x}^{\sharp}$
in the equality (\ref{eq:Dy}). Then we have $\beta(\widetilde{x},y)=\beta(\widetilde{x},\widetilde{x}^{\sharp})=3N_{x}(\widetilde{x})$
by (\ref{eq:normxy}). By (\ref{eq:xxshay}), we also have 
\[
\widetilde{x}\bullet y\bullet\widetilde{y}=\widetilde{x}\bullet\widetilde{x}^{\sharp}\bullet\widetilde{y}=1/6(-\beta(\widetilde{x},\widetilde{y})\widetilde{x}^{\sharp}+N_{x}(\widetilde{x})\widetilde{y}).
\]
Therefore the equality (\ref{eq:Dy}) with $y=\widetilde{x}^{\sharp}$
becomes 
\[
1/6\,\beta(\widetilde{x},\widetilde{y})\widetilde{x}^{\sharp}-3/2\,N_{x}(\widetilde{x})\widetilde{y}-1/2\,(-\beta(\widetilde{x},\widetilde{y})\widetilde{x}^{\sharp}+N_{x}(\widetilde{x})\widetilde{y})=0,
\]
equivalently, we have $\widetilde{y}=\frac{\beta(\widetilde{x},\widetilde{y})}{3N_{x}(\widetilde{x})}\widetilde{x}^{\sharp}$.
Thus (2) follows. 

Now we show (2) implies (1). Assume that (2) holds. We write $\widetilde{y}=\alpha\widetilde{x}^{\sharp}$
for some $\alpha\in\mathsf{k}$ (actually, computing $\beta(\widetilde{x},\widetilde{y}),$
we have $\alpha=\frac{\beta(\widetilde{x},\widetilde{y})}{3N_{x}(\widetilde{x})}$).
Inserting $\widetilde{y}=\alpha\widetilde{x}^{\sharp}$ in the l.h.s.~of
(\ref{eq:Dy}), we see that the equality (\ref{eq:Dy}) holds for
any $y\in V_{y}$ by Lemma \ref{lem:xshaId}.
\end{proof}
We recall that we set $n:=\dim V_{x}=\dim V_{y}$. By non-degeneracy
of $\beta$, we immediately obtain the following from Lemma \ref{lem:n-1dim}:
\begin{cor}
\label{cor:n-1} For $\widetilde{x}\in V_{x}$ with $N_{x}(\widetilde{x})\not=0$,
the entries of the bi-linear map 
\[
\Delta_{x}(*,*,\widetilde{x})\colon V_{x}\times V_{y}\to V_{x}
\]
 generate an $(n-1)$-dimensional vector space. 
\begin{example}
\label{exa:twoJ}As mentioned in the beginning of Subsection \ref{subsec:Jordan-algebraic-description},
a basic example of an FTS is the one containing two copies of a Jordan
algebra of a cubic form. For such an example, $V_{x}=V_{y}$ and this
has the structure of the Jordan algebra of the cubic form compatible
with the quartic form $F$ and the skew-symmetric form $\omega$.
For the Jordan algebra $V:=V_{x}=V_{y}$, the cubic form $N:=N_{x}=N_{y}$
is the associated cubic form and the bilinear form $\beta$ is the
associated bi-linear trace (see \cite[p.314]{Fe} for more details).

\begin{enumerate}

\item Let $V$ be a $3$-dimensional vector space with coordinates
$x_{1},x_{2},x_{3}$. The cubic form $x_{1}x_{2}x_{3}$ define a Jordan
algebra structure on $V$, whose $\sharp$-mapping is $V\ni(x_{1},x_{2},x_{3})\mapsto(x_{2}x_{3},x_{1}x_{3},x_{1}x_{2})\in V,$
and the bi-linear trace is $\beta(x,y)=x_{1}y_{1}+x_{2}y_{2}+x_{3}y_{3}$
for $x=(x_{1},x_{2},x_{3})$ and $y=(y_{1},y_{2},y_{3}).$ By an explicit
calculation, we see that the condition (\ref{eq:forany2}) is reduced
to the two equations: $\,2x_{1}y_{1}-x_{2}y_{2}-x_{3}y_{3}=0,x_{1}y_{1}-2x_{2}y_{2}+x_{3}y_{3}=0.$
Then we can write down the 9 equations of $\mathfrak{R}$ as follows:
\begin{align*}
 & sx_{1}=y_{2}y_{3},\,sx_{2}=y_{1}y_{3},\,sx_{3}=y_{1}y_{2},\\
 & ty_{1}=x_{2}x_{3},\,ty_{2}=x_{1}x_{3},\,ty_{3}=x_{1}x_{2},\\
 & st=1/3(x_{1}y_{1}+x_{2}y_{2}+x_{3}y_{3}),\\
 & 2x_{1}y_{1}-x_{2}y_{2}-x_{3}y_{3}=0,\,x_{1}y_{1}-2x_{2}y_{2}+x_{3}y_{3}=0.
\end{align*}
Tidying up these equations, we see that these are derived from a $2\times2\times2$
hypermatrix and then $\mathfrak{R}$ is the affine cone over the Segre
embedded $\mP^{1}\times\mP^{1}\times\mP^{1}$.

\item The space of $3\times3$ matrices $M(3,\mathsf{k})$ is a Jordan
algebra of a cubic form. The Jordan product $\cdot_{J}$ is defined
as follows: 
\[
X\cdot_{J}Y:=1/2(XY+YX)\ \text{for\ }X,Y\in M(3,\mathsf{k}).
\]
The associated cubic form $N$ is the determinant of a matrix, and
for $X\in M(3,\mathsf{k}),$ $X^{\sharp}$ is the adjoint matrix.
The bi-linear trace $\beta$ is defined as follows:
\[
\beta(X,Y)=\text{tr}(XY)\ \text{for\ }X,Y\in M(3,\mathsf{k}).
\]
In this setting, it is known that $\mathfrak{R}$ is the Grassmannian
${\rm G}(3,6)$ (cf.~\cite{LM}) and we can compute the defining
equation $\mathfrak{R}$ as follows:
\[
sX=Y^{\sharp},\,tY=X^{\sharp},\,XY=YX,\,XY=stI,
\]
where $I$ is the $3\times3$ identity matrix.

\item The space of $3\times3$ symmetric matrices ${\rm Sym}(3,\mathsf{k})$
is also a Jordan algebra of a cubic form, and is a Jordan subalgebra
of $M(3,\mathsf{k})$. In this setting, $\mathfrak{R}$ is the $6$-dimensional
symplectic Grassmannian ${\rm Sp}(3,6)$ and we can compute the defining
equation of $\mathfrak{R}$ as follows:
\[
sX=Y^{\sharp},\,tY=X^{\sharp},\,XY=stI
\]
(cf.\cite[p.32]{I}).

\end{enumerate}
\end{example}

\end{cor}

\subsection{More on FTS and the scheme $\mathfrak{R}$}

Finally in this section, we add further properties of FTS and the
scheme $\mathfrak{R}$.

Explicit descriptions of the quartic form\textit{ $\widetilde{F}$}
and the symplectic form $\omega$ are given as follows (cf.\cite[(21) in Sec.7, p.118--119]{Cl},\cite[p.321]{Fe},
\cite[Prop.5.5]{LM}). The proof is same as that given in \cite[p.321]{Fe}
with the above formulation. 
\begin{prop}
\label{prop:Explicit}For $p=(s,t,x,y)$ and $p'=(s',t',x',y')\in V$,
it holds that
\begin{align*}
 & \widetilde{F}(p)=8\left(\beta\left(y^{\sharp},x^{\sharp}\right)-sN_{x}(x)-tN_{y}(y)\right)-2\left(st-\beta(x,y)\right)^{2},\\
 & \omega(p,p')=\beta(x,y')-\beta(x',y)+st'-s't.
\end{align*}
\end{prop}

From this, we see that $\omega$ is nondegenerate if and only if $\beta$
is so. We use this proposition to show Proposition \ref{prop:converse}
below. 

\vspace{3pt}

We see that $\mathfrak{R}$ contains an open subset isomorphic to\textit{
$\mathsf{k}^{*}\times\mA^{n}$} as follows (cf.\cite[p.114]{Cl},
\cite[(4.3)]{Mu2}):
\begin{cor}
\label{cor:affinespace}The following two assertions hold: 

\begin{enumerate}

\item An element $p=(s,t,x,y)$ with $s\not=0$ is strictly regular
if and only if it holds that $x=s^{-1}y^{\sharp},t=(3s)^{-1}\beta(x,y)$.

\item An element $p=(s,t,x,y)$ with $t\not=0$ is strictly regular
if and only if it holds that $y=t^{-1}x^{\sharp},s=(3t)^{-1}\beta(x,y)$.

\end{enumerate}

In particular, $\mathfrak{R}\cap\{s\not=0\}$ and $\mathfrak{R}\cap\{t\not=0\}$
are isomorphic to $\mathsf{k}^{*}\times\mA^{n}$.
\end{cor}

\begin{proof}
Since we can prove these assertions similarly, we only show (1). The
only if part follows immediately by Corollary \ref{cor:streg}. To
show the if part, we assume that $x=s^{-1}y^{\sharp},t=(3s)^{-1}\beta(x,y)$
for $p=(s,t,x,y)$. By Corollary \ref{cor:streg}, it suffices to
check that 
\begin{equation}
\left((3s)^{-1}\beta(s^{-1}y^{\sharp},y)\right)y-(s^{-1}y^{\sharp})^{\sharp}=0\label{eq:check1}
\end{equation}
and 
\begin{equation}
\Delta_{x}(s^{-1}y^{\sharp},y,\widetilde{x})=\Delta_{y}(s^{-1}y^{\sharp},y,\widetilde{y})=0\,\,\text{for any }\widetilde{p}=(\widetilde{s},\widetilde{t,}\widetilde{x},\widetilde{y}).\label{eq:check2}
\end{equation}
The equality (\ref{eq:check1}) follows from Propositions \ref{prop:NxNy}
and \ref{prop:=000023=000023}, and the equality (\ref{eq:check2})
follows from Lemma \ref{lem:xshaId}.
\end{proof}

\subsection{A coordinatization of $\mathfrak{R}$ with parameters\label{subsec:n=00003D3}}

In this subsection, we further assume that 
\begin{align*}
\text{\ensuremath{}} & n=\dim V_{x}=\dim V_{y}=3.
\end{align*}
We show that the FTS $V$ has a coordinatization with parameters as
in the following theorem, which can be seen as a generalization of
Example \ref{exa:twoJ} (1):
\begin{thm}
\label{thm:EqF22} Assume that $N_{x}$ is not identically zero on
$V_{x}$. The equations (\ref{eq:forany2}) as in Theorem \ref{thm:stregoneeq}
are reduced to the two equations 
\[
{\empty}^{t}\!xPy=0,{\empty}^{t}\!xQy=0,
\]
where $P,Q$ are $3\times3$ $k$-matrices, and $x,y$ are considered
as column vectors ${\empty}^{t}\!(x_{1},x_{2},x_{3})$ and ${\empty}^{t}\!(y_{1},y_{2},y_{3})$.
Moreover, there exist nonzero constants $\mu,\nu\in\mathsf{k}$ such
that 
\begin{equation}
x^{\sharp}=\mu({\empty}^{t}\!Px\times{\empty}^{t}\!Qx),\,y^{\sharp}=\nu(Py\times Qy),\label{eq:xysharp}
\end{equation}
 and, by replacing $s,t$ with $\mu s,\nu$t respectively, the equation
of $\mathfrak{R}$ is reduced to 
\begin{align*}
 & sx=Py\times Qy,\,ty={\empty}^{t}\!Px\times{\empty}^{t}Qx,\\
 & st=\frac{1}{3\mu\nu}\beta(x,y),\\
 & {\empty}^{t}\!xPy=0,{\empty}^{t}\!xQy=0,
\end{align*}
where ${\empty}^{t}\!Px\times{\empty}^{t}\!Qx$ is the cross product
of the two column vectors ${\empty}^{t}\!Px$ and ${\empty}^{t}\!Qx$,
and $Py\times Qy$ is similarly defined. In particular, the scheme
$\mathfrak{R}$ in this case is defined by $9$ quadratic forms.
\end{thm}

\begin{proof}
Since $N_{x}$ is not identically zero on $V_{x}$, we may take $\widetilde{x}\in V_{x}$
with $N_{x}(\widetilde{x})\not=0$. Since $\dim V_{x}=\dim V_{y}=3$,
the entries of the bi-linear map 
\[
\Delta_{x}(*,*,\widetilde{x})\colon V_{x}\times V_{y}\to V_{x}
\]
generate a $2$-dimensional vector space by Corollary \ref{cor:n-1}.
We denote a basis of this vector space by ${\empty}^{t}\!xPy$ and
${\empty}^{t}\!xQy$ with some $3\times3$ matrices $P$ and $Q$
with entries in $\mathsf{k}$. By Lemma \ref{lem:xshaId}, we have
\begin{align}
{\empty}^{t}\!xPx^{\sharp} & =0,{\empty}^{t}\!xQx^{\sharp}=0\,\text{for any}\,x\in V_{x},\label{eq:xP}\\
{\empty}^{t}\!y^{\sharp}Py & =0,{\empty}^{t}\!y^{\sharp}Qy=0\,\text{for any}\,y\in V_{y}.\label{eq:yP}
\end{align}
By (\ref{eq:xP}), there exists a rational function $\mu(x)$ such
that 
\begin{equation}
x^{\sharp}=\mu(x)({\empty}^{t}\!Px\times{\empty}^{t}\!Qx).\label{eq:mux}
\end{equation}
 Since all the entries of both $x^{\sharp}$ and ${\empty}^{t}\!Px\times{\empty}^{t}\!Qx$
are quadratic forms, the degree of $\mu(x)$ is 0. By Corollary \ref{cor:nocommonfac=000023},
$\mu(x)$ is a nonzero constant, hence we will denote this by $\mu$.
Similarly, we can show that there exists a nonzero constant $\nu\in k$
such that
\[
y^{\sharp}=\nu(Py\times Qy).
\]
\end{proof}
Under the situation of Theorem \ref{thm:EqF22}, we can write down
several data of FTS using the entries of $x,y,P,Q$. Here we treat
some of them as follows:

\vspace{5pt}

\noindent \textbf{Recipe to write down $N_{x}(x),N_{y}(y),$ and
$\beta(x,y)$. }

\begin{enumerate}[(1)]

\item By (\ref{eq:xysharp}), we have 
\begin{align*}
 & (x^{\sharp})^{\sharp}=\mu^{2}\nu\left(P({\empty}^{t}\!Px\times{\empty}^{t}\!Qx)\times Q({\empty}^{t}\!Px\times{\empty}^{t}\!Qx)\right),\\
 & (y^{\sharp})^{\sharp}=\mu\nu^{2}\left({\empty}^{t}\!P(Py\times Qy)\times{\empty}^{t}\!Q(Py\times Qy)\right).
\end{align*}
From these, we may write down the cubic forms $N_{x}(x)$ and $N_{y}(y)$
by Proposition \ref{prop:=000023=000023}. 

\item By (1) and Corollary \ref{cor:partialN}, we may write down
the bi-linear trace $\beta(x,y)$ since $N_{x}(x)\not=0$ for some
$x\in V_{x}$ by the assumption, and then the subset $\{x^{\sharp}\mid x\in V_{x}\}\subset V_{y}$
is dense in $V_{y}$ by Proposition \ref{prop:=000023=000023}.

\end{enumerate}

\vspace{5pt}

Suggested by this recipe, we have a converse statement to Theorem
\ref{thm:EqF22}. The verification of the following lemma and proposition
is straightforward.
\begin{lem}
\label{lem:Converse}Let $x$, $y$ be $3$-dimensional column vectors
and $P$, $Q$ $3\times3$ matrices. We consider all the entries of
$x$, $y$, $P$, $Q$ as variables of a $\mathsf{k}$-polynomial
ring. The following assertions hold:

\begin{enumerate}[(1)]

\item There exist the uniquely determined cubic forms $N_{x}^{PQ}(x)$
and $N_{y}^{PQ}(y)$ such that 
\begin{align*}
 & \left(P({\empty}^{t}\!Px\times{\empty}^{t}\!Qx)\times Q({\empty}^{t}\!Px\times{\empty}^{t}\!Qx)\right)=N_{x}^{PQ}(x)x,\\
 & \left({\empty}^{t}\!P(Py\times Qy)\times{\empty}^{t}\!Q(Py\times Qy)\right)=N_{y}^{PQ}(y)y.
\end{align*}
\item Let $x'$, $y'$ be $3$-dimensional column vectors. There
exists the uniquely determined bi-linear trace $\beta^{PQ}(x,y)$
such that 
\begin{equation}
\partial_{x'}N_{x}^{PQ}(x)=\beta^{PQ}\left(x',{\empty}^{t}\!Px\times{\empty}^{t}\!Qx\right),\partial_{y'}N_{y}^{PQ}(y)=\beta^{PQ}\left(Py\times Qy,y'\right).\label{eq:normxypartial-1}
\end{equation}

\end{enumerate}
\end{lem}

\begin{prop}
\label{prop:converse}Let $V_{x}$ and $V_{y}$ be $3$-dimensional
$k$-vector spaces. We write elements of $V_{x}$ and $V_{y}$ as
column vectors. Let $P,Q$ be $3\times3$ $\mathsf{k}$-matrices.
We define 
\[
x^{\sharp}:={\empty}^{t}\!Px\times{\empty}^{t}\!Qx,\,y^{\sharp}:=Py\times Qy\ \text{for}\ x\in V_{x},\,y\in V_{y},
\]
and the $\sharp$-product by (\ref{eq:bisha}). Let the cubic forms
$N_{x}^{PQ}(x)$, $N_{y}^{PQ}(y)$ and the bi-linear form $\beta^{PQ}(x,y)$
($x\in V_{x}$, $y\in V_{y}$) be as defined in Lemma \ref{lem:Converse}.
Let $V:=\mathsf{k}\oplus\mathsf{k}\oplus V_{x}\oplus V_{y}$. We further
define a skew bi-linear form $\omega$ on $V$ as follows: for $p=(s,t,x,y)$
and $p'=(s',t',x',y')\in V$,
\begin{align*}
\\
 & \omega(p,p'):=\beta^{PQ}(x,y')-\beta^{PQ}(x',y)+st'-s't.
\end{align*}
 Finally, we define a tri-linear product on $V$ by completely linearizing
the formula (\ref{eq:ppp}) with $x^{\sharp}$, $y^{\sharp}$, $N_{x}^{PQ}(x)$,
$N_{y}^{PQ}(y)$ and $\beta^{PQ}(x,y)$. The tri-linear product satisfies
the axioms as in Definition \ref{def:Triple} with respect to $\omega$,
and the symmetric $4$-linear form as in $(A2)$ is the complete linearization
of the following quartic form $\widetilde{F}^{PQ}:$ 
\[
\widetilde{F}^{PQ}(p):=8\left(\beta^{PQ}\left(y^{\sharp},x^{\sharp}\right)-sN_{x}^{PQ}(x)-tN_{y}^{PQ}(y)\right)-2\left(st-\beta^{PQ}(x,y)\right)^{2}.
\]
In particular the vector space $V$ is the FTS with respect to $\omega$
and this tri-linear product.
\end{prop}

\begin{rem}
\label{rem:Dbeta}Let $D_{\beta}$ be the determinant of the matrix
defining the bi-linear form $\beta^{PQ}(x,y)$ as in Proposition \ref{prop:converse};
$\beta^{PQ}(x,y)$ is nondegenerate if and only if $D_{\beta}\not=0$.
\end{rem}

\begin{cor}
\label{cor:P1P1P1}Under the situation as in Proposition \ref{prop:converse},
assume that the base field $\mathsf{k}$ is algebraically closed and
$\beta^{PQ}(x,y)$ is nondegenerate, namely, $D_{\beta}\not=0$. Let
$\mathfrak{R}^{PQ}$ be the affine scheme defined by 
\begin{align}
 & sx=Py\times Qy,\,ty={\empty}^{t}\!Px\times{\empty}^{t}Qx,\label{eq:1st22}\\
 & st=1/3\,\beta^{PQ}(x,y),\label{eq:2nd22}\\
 & {\empty}^{t}\!xPy=0,{\empty}^{t}\!xQy=0.\label{eq:3rd22}
\end{align}
The affine scheme $\mathfrak{R}^{PQ}$ is isomorphic to the affine
cone over the Segre embedded $\mP^{1}\times\mP^{1}\times\mP^{1}.$
\end{cor}

\begin{proof}
Since $\beta^{PQ}(x,y)$ is nondegenerate, $V$ is simple as an FTS
by \cite[Thm.2.1]{Fe}. Since $(1,0,0,0)$, $(0,1,0,0)\in V$ are
supplementary strictly regular elements, $V$ is reduced. Therefore,
by \cite[p.314 and Thm.5.1]{Fe}, $V$ is constructed from a Jordan
algebra $J$ as in (i) or (ii) of \cite[Thm.5.1]{Fe}. 

We show (i) holds. Assume for a contradiction that (ii) holds. Then
\begin{equation}
{\empty}^{t}\!Px\times{\empty}^{t}Qx=0\ \text{for any}\ x\in V_{x}\label{eq:PxQx}
\end{equation}
 by \cite[p.314]{Fe}. If $P=0$ or $Q=0$, then $\beta^{PQ}(x,y)$
is identically zero, a contradiction. Therefore, changing the coordinates
of $V_{x}$ and $V_{y}$, we may assume that $P={\rm diag}\,(1,1,1),{\rm diag}\,(1,1,0),\ \text{or}\ {\rm diag\,}(1,0,0)$,
where ${\rm diag}\,(a,b,c)$ is the diagonal matrix with $a,b,c$
as the $(1,1)$, $(2,2)$, $(3,3)$-entries respectively. Then, from
(\ref{eq:PxQx}), we see that $Q=\alpha P$ for some $\alpha\in\mathsf{k}$
in the first or the second case, or $Q=\left(\begin{array}{ccc}
* & 0 & 0\\
* & 0 & 0\\
* & 0 & 0
\end{array}\right)$ in the third case. In any case, we see that $\beta^{PQ}(x,y)$ is
identically zero, a contradiction. Therefore (i) must hold. 

In this case, nondegeneracy of $\beta^{PQ}(x,y)$ implies that the
Jordan algebra $J$ does not contain an absolute zero divisor in the
sense of \cite[p.93]{Ra} (\cite[Def.5.3.2, Prop.5.3.3, Ex.5.3.6]{Mc2}).
Therefore, (i), (ii) or (iii) of \cite[Thm.1]{Ra} holds. The Jordan
algebra $J$ is not a division algebra since $N_{x}=0$ for some nonzero
$x\in V_{x}$. Thus (ii) or (iii) holds. Then, by the proof of ibid.,
the Jordan algebra $J$ has the Peirce decomposition $J=J_{1}\oplus J_{1/2}\oplus J_{0}$
with respect to a primitive idempotent. We have $\dim J_{1}=1$ and
$\dim J_{0}\geq1$. If (iii) holds, then $\dim J_{1/2}\geq2$ by \cite[p.98, 5 and 6th line from the bottom]{Ra},
hence $\dim V_{x}=\dim J\geq4$, a contradiction. Therefore (ii) holds,
namely, $J=\mathsf{k}\oplus J(q)$, where $J(q)$ is a Jordan algebra
of a quadratic form $q$. By \cite[p.506, Ex.2]{Mc1}, $q$ is nondegenerate
since so is $\beta^{PQ}(x,y)$. Therefore, by ibid. and the assumption
that $\mathsf{k}$ is algebraically closed, we see that $N_{x}$ is
a product of three linearly independent linear forms. Now, by Example
\ref{exa:twoJ} (1), $\mathfrak{R}^{PQ}$ is isomorphic to the affine
cone over $\mP^{1}\times\mP^{1}\times\mP^{1}.$
\end{proof}

\part{Key varieties for $\mQ$-Fano threefolds\protect 
}\textit{Hereafter we work over $\mC$, the complex number field, throughout
the paper.}

\section{\textbf{The affine scheme $\mathfrak{F}_{\mA}^{22}$ \label{sec:Affine-variety-22}}}
\begin{defn}
\label{def:In-the-affineF22}In the affine $26$-space whose coordinates
are $s,t,$ and the entries of $x,y,P,Q$, we define $\mathfrak{F}_{\mathbb{A}}^{22}$
to be the scheme with the equations (\ref{eq:1st22}), (\ref{eq:2nd22})
and (\ref{eq:3rd22}). We call the entries of $P,Q$ \textit{the parametric
coordinates.}
\end{defn}

\begin{rem}
The equations of $\mathfrak{F}_{\mathbb{A}}^{22}$ is a specialization
of those given in \cite[Subsec.5.7]{P2}. These are also derived in
\cite{NP}. The main advantage here is that the meaning of the equations
is quite clear in view of the theory of FTS; especially, the equation
(\ref{eq:2nd22}) is too complicated to write down fully but we can
write it conceptually as above. 
\end{rem}

In the following Sections \ref{sec:KeyEquation14}--\ref{sec:The--Cluster-variety},
setting some of the parametric coordinates of $\mathfrak{F}_{\mathbb{A}}^{22}$
to constants, or imposing linear relations on them, we will obtain
several affine varieties whose weighted projectivizations produce
examples of prime $\mQ$-Fano 3-folds of anti-canonical codimension
4. We refer to Section \ref{sec:Affine-variety-Papadakis} for another
aspect of $\mathfrak{F}_{\mathbb{A}}^{22}$.

We set 
\[
M_{x}:=\left(\begin{array}{c}
{\empty}^{t}\!xP\\
{\empty}^{t}\!xQ
\end{array}\right),M_{y}:=\left(\begin{array}{c}
{\empty}^{t}\!y\empty^{t}\!P\\
{\empty}^{t}\!y\empty^{t}\!Q
\end{array}\right).
\]

\section{\textbf{A specialization of $\mathfrak{F}_{\mathbb{A}}^{22}$--an
affine variety $\mathfrak{U}_{\mathbb{A}}^{14}$ with an ${\rm SL}_{2}\times{\rm SL}_{2}$-action}--
\label{sec:KeyEquation14}}

In this section, we define a closed subscheme $\mathfrak{U}_{\mathbb{A}}^{14}$
of the affine scheme $\mathfrak{F}_{\mathbb{A}}^{22}$ and investigate
its properties. We show that a weighted projectivization of $\mathfrak{U}_{\mathbb{A}}^{14}$
give an example of a key variety for a prime $\mQ$-Fano 3-fold of
anticanonical codimension 4 which is not obtained from the $G_{2}^{(4)}$-cluster
variety. 

\subsection{Definition\label{subsec:DefinitionU}}
\begin{defn}
We set 
\[
\mathfrak{U}_{\mathbb{A}}^{14}:=\mathfrak{F}_{\mathbb{A}}^{22}\cap\{p_{13}=q_{23}=1,p_{23}=p_{33}=q_{13}=q_{33}=0,p_{11}=-q_{21},p_{12}=-q_{22}\}.
\]
\end{defn}

We use the following notation for the entries of the matrices $P$
and $Q$, by which a group action on $\mathfrak{U}_{\mathbb{A}}^{14}$
will be more visible (see Subsection \ref{subsec:SL2SL2}): 

\[
P=\left(\begin{array}{ccc}
a_{11} & b_{11} & 1\\
a_{12} & b_{12} & 0\\
c_{11} & c_{12} & 0
\end{array}\right),Q=\left(\begin{array}{ccc}
a_{21} & b_{21} & 0\\
-a_{11} & -b_{11} & 1\\
c_{21} & c_{22} & 0
\end{array}\right).
\]
Setting further
\begin{align*}
 & A:=\left(\begin{array}{cc}
a_{11} & a_{12}\\
a_{21} & -a_{11}
\end{array}\right),B:=\left(\begin{array}{cc}
b_{11} & b_{12}\\
b_{21} & -b_{11}
\end{array}\right),C:=\left(\begin{array}{cc}
c_{11} & c_{12}\\
c_{21} & c_{22}
\end{array}\right),\\
 & \widehat{x}=\left(\begin{array}{c}
x_{1}\\
x_{2}
\end{array}\right),\widehat{y}=\left(\begin{array}{c}
y_{1}\\
y_{2}
\end{array}\right),
\end{align*}
we can write 

\begin{align*}
M_{x} & =\left(\begin{array}{c}
{\empty}^{t}\!xP\\
{\empty}^{t}\!xQ
\end{array}\right)=\left(\begin{array}{ccc}
A\widehat{x} & B\widehat{x} & \widehat{x}\end{array}\right)+x_{3}\left(\begin{array}{cc}
C & 0\end{array}\right).
\end{align*}
Explicitly,
\[
M_{x}=\left(\begin{array}{ccc}
a_{11}x_{1}+a_{12}x_{2}+c_{11}x_{3} & b_{11}x_{1}+b_{12}x_{2}+c_{12}x_{3} & x_{1}\\
a_{21}x_{1}-a_{11}x_{2}+c_{21}x_{3} & b_{21}x_{1}-b_{11}x_{2}+c_{22}x_{3} & x_{2}
\end{array}\right).
\]
We have also

\[
M_{y}=\left(\begin{array}{c}
{\empty}^{t}\!y\empty^{t}\!P\\
{\empty}^{t}\!y\empty^{t}\!Q
\end{array}\right)=\left(\begin{array}{ccc}
a_{11}y_{1}+b_{11}y_{2}+y_{3} & a_{12}y_{1}+b_{12}y_{2} & c_{11}y_{1}+c_{12}y_{2}\\
a_{21}y_{1}+b_{21}y_{2} & -a_{11}y_{1}-b_{11}y_{2}+y_{3} & c_{21}y_{1}+c_{22}y_{2}
\end{array}\right).
\]
We denote by $\mA_{\mathfrak{U}}^{18}$ the affine space whose coordinates
are $s,t$ and the entries of $A,B,C,x,y$. 

\subsection{${\rm SL}_{2}\times{\rm SL}_{2}$-action\label{subsec:SL2SL2}}

We see that $\mathfrak{U}_{\mathbb{A}}^{14}$ admits the following
group action, which will be helpful to investigate properties of $\mathfrak{U}_{\mathbb{A}}^{14}$:
\begin{prop}
The scheme $\mathfrak{U}_{\mathbb{A}}^{14}$ is preserved by the following
group actions on the affine space $\mA_{\mathfrak{U}}^{18}$ of the
two groups $({\rm SL_{2}})^{x}$, $({\rm SL_{2}})^{y}$ isomorphic
to ${\rm SL_{2}}$, and they define an action of the group ${\rm SL}_{2}\times{\rm SL}_{2}$
on $\mathfrak{U}_{\mathbb{A}}^{14}$:
\begin{align*}
({\rm SL_{2}})^{x}:\text{For \ensuremath{g}\ensuremath{\in({\rm SL}_{2})^{x}}, } & A\mapsto gAg^{-1},\,B\mapsto gBg^{-1},\widehat{x}\mapsto g\widehat{x},\,C\mapsto gC,\\
 & x_{3}\mapsto x_{3},\,y\mapsto y,\,s\mapsto s,\,t\mapsto t,\\
({\rm SL_{2}})^{y}:\text{For \ensuremath{h\in({\rm SL}_{2})^{y}},} & \left(\begin{array}{cc}
a_{11} & b_{11}\\
a_{21} & b_{21}
\end{array}\right)\mapsto\left(\begin{array}{cc}
a_{11} & b_{11}\\
a_{21} & b_{21}
\end{array}\right)h^{-1},\\
 & \left(\begin{array}{cc}
a_{12} & b_{12}\\
-a_{11} & -b_{11}
\end{array}\right)\mapsto\left(\begin{array}{cc}
a_{12} & b_{12}\\
-a_{11} & -b_{11}
\end{array}\right)h^{-1},\\
 & \widehat{y}\mapsto h\widehat{y},\,C\mapsto Ch^{-1},\\
 & y_{3}\mapsto y_{3},\,x\mapsto x,\,s\mapsto s,\,t\mapsto t.
\end{align*}
\end{prop}

\begin{proof}
Note that $M_{x}y={\empty}^{t}\!(M_{y}x)=\left(\begin{array}{c}
{\empty}^{t}\!xPy\\
{\empty}^{t}\!xQy
\end{array}\right).$ By straightforward calculations, we see the following:

\begin{itemize}

\item $M_{x}$ is mapped to $gM_{x}$ by $\ensuremath{g}\ensuremath{\in({\rm SL}_{2}})^{x}$,
and $M_{y}$ is invariant for the action of $({\rm SL}_{2})^{y}$. 

\item $x^{\sharp}$ is invariant for the action of $({\rm SL}_{2})^{x}$
and is equivariant to $y$ for the action of $({\rm SL}_{2})^{y}$.
$y^{\sharp}$ is equivariant to $x$ for the action of $({\rm SL}_{2})^{x}$
and is invariant for the action of $({\rm SL}_{2})^{y}$.

\item $\beta(x,y)$ is invariant for the actions of $({\rm SL}_{2})^{x}$
and $({\rm SL}_{2})^{y}$. 

\item The actions of $({\rm SL}_{2})^{x}$ and $({\rm SL}_{2})^{y}$
are commutative.

\end{itemize}

Therefore we have the group action on $\mathfrak{U}_{\mathbb{A}}^{14}$
as in the statement.
\end{proof}

\subsection{Weights on the variables and the equations \label{subsec:posgr}}

We assign weights on the variables of the polynomial ring $S_{\mathfrak{U}}$
such that all the $9$ equations of $\mathfrak{U}_{\mA}^{14}$ are
homogeneous. Moreover, we assume that all the variables are not zero
allowing some of them to be constants. Then it is easy to derive the
following relations between the weights on the variables of $S_{\mathfrak{U}}$:

\begin{align*}
w(a_{11}) & =-w(y_{1})+w(y_{3}),\,w(a_{12})=w(x_{1}y_{3})-w(x_{2}y_{1}),\,w(a_{21})=w(x_{2}y_{3})-w(x_{1}y_{1}),\\
w(b_{11}) & =-w(y_{2})+w(y_{3}),\,w(b_{12})=w(x_{1}y_{3})-w(x_{2}y_{2}),\,w(b_{21})=w(x_{2}y_{3})-w(x_{1}y_{2}),\\
w(c_{11}) & =w(x_{1}y_{3})-w(x_{3}y_{1}),\,w(c_{12})=w(x_{1}y_{3})-w(x_{3}y_{2}),\\
w(c_{21}) & =w(x_{2}y_{3})-w(x_{3}y_{1}),\,w(c_{22})=w(x_{2}y_{3})-w(x_{3}y_{2}),\\
w(s) & =-w(x_{3})+2w(y_{3}),\,w(t)=w(x_{1}x_{2}y_{3})-w(y_{1}y_{2}).
\end{align*}

\begin{example}
\label{exa:In-Subsection-,WT U}In Subsection \ref{subsec:An-example-of fano 1},
we use the following weights on the coordinates:
\begin{align*}
 & w(A)=w(B)=w(C)=\left(\begin{array}{cc}
1 & 1\\
1 & 1
\end{array}\right),\\
 & w(x)=\left(\begin{array}{c}
1\\
1\\
1
\end{array}\right),w(y)=\left(\begin{array}{c}
1\\
1\\
2
\end{array}\right),\\
 & w(s)=3,w(t)=2.
\end{align*}
Let 
\begin{align*}
 & \mathfrak{U}_{\mP}^{13}\subset\mP(1^{15},2^{2},3)
\end{align*}
be the weighted projectivization of $\mathfrak{U}_{\mathbb{A}}^{14}$
by these weights on the coordinates. In Subsection \ref{subsec:An-example-of fano 1},
we show that $\mathfrak{U}_{\mP}^{13}$ is a key variety of a prime
$\mQ$-Fano 3-fold of anticanonical codimension 4 with No.20544. Existence
of a positive grading of $\mathfrak{U}_{\mathbb{A}}^{14}$ also plays
a role in the proof of Proposition \ref{prop:Jerry14}. 
\end{example}

\subsection{Charts and singular locus\label{subsec:Descriptrion-of-charts}}

For a coordinate $*$, we call the open subset of $\mathfrak{U}_{\mA}^{14}$
with $*\not=0$ \textit{the $*$-chart}. We describe the $*$-chart
such that $*$ is one of the entries of $x,y$, or $s$, $t$.

\vspace{3pt}

\noindent \underline{\bf $x_1$-chart}: We note that $\mathfrak{U}_{\mA}^{14}\cap\{x_{1}\not=0\}$
is isomorphic to $(\mathfrak{U}_{\mA}^{14}\cap\{x_{1}=1\})\times\mA^{1*}$
by the map $(x,y,s,t,A,B,C)\mapsto\left((x_{1}^{-1}x,x_{1}^{-1}y,x_{1}^{-1}s,x_{1}^{-1}t,A,B,C),x_{1}\right)$.
This is because all the equations of $\mathfrak{U}_{\mA}^{14}$ are
quadratic when we consider the entries of $A,B,C$ are constants.
Therefore it suffices to describe $\mathfrak{U}_{\mA}^{14}\cap\{x_{1}=1\}$.
Solving regularly the $9$ equations of $\mathfrak{U}_{\mA}^{14}$
setting $x_{1}=1$, we see that the $9$ equations are reduced to
the following 4 equations: 

\begin{align*}
y_{3}= & -a_{11}y_{1}-a_{12}x_{2}y_{1}-c_{11}x_{3}y_{1}-b_{11}y_{2}-b_{12}x_{2}y_{2}-c_{12}x_{3}y_{2},\\
s= & 2a_{11}c_{11}y{}_{1}^{2}+a_{12}c_{21}y{}_{1}^{2}+a_{12}c_{11}x_{2}y{}_{1}^{2}+c_{11}^{2}x_{3}y{}_{1}^{2}+2b_{11}c_{11}y_{1}y_{2}+\\
 & 2a_{11}c_{12}y_{1}y_{2}+b_{12}c_{21}y_{1}y_{2}+a_{12}c_{22}y_{1}y_{2}+b_{12}c_{11}x_{2}y_{1}y_{2}+a_{12}c_{12}x_{2}y_{1}y_{2}+\\
 & 2c_{11}c_{12}x_{3}y_{1}y_{2}+2b_{11}c_{12}y{}_{2}^{2}+b_{12}c_{22}y{}_{2}^{2}+b_{12}c_{12}x_{2}y{}_{2}^{2}+c_{12}^{2}x_{3}y{}_{2}^{2},\\
b_{21}= & 2b_{11}x_{2}+b_{12}x{}_{2}^{2}-c_{22}x_{3}+c_{12}x_{2}x_{3}-ty_{1,}\\
a_{21}= & 2a_{11}x_{2}+a_{12}x{}_{2}^{2}-c_{21}x_{3}+c_{11}x_{2}x_{3}+ty_{2}.
\end{align*}

\vspace{3pt}

\noindent \underline{\bf $x_2$-chart}: Similarly to the $x_{1}$-chart,
we have only to describe $\mathfrak{U}_{\mA}^{14}\cap\{x_{2}=1\}$
as follows:

\begin{align*}
y_{3}= & a_{11}y_{1}-a_{21}x_{1}y_{1}-c_{21}x_{3}y_{1}+b_{11}y_{2}-b_{21}x_{1}y_{2}-c_{22}x_{3}y_{2},\\
s= & a_{21}c_{11}y{}_{1}^{2}-2a_{11}c_{21}y{}_{1}^{2}+a_{21}c_{21}x_{1}y{}_{1}^{2}+c_{21}^{2}x_{3}y{}_{1}^{2}+b_{21}c_{11}y_{1}y_{2}+\\
 & a_{21}c_{12}y_{1}y_{2}-2b_{11}c_{21}y_{1}y_{2}-2a_{11}c_{22}y_{1}y_{2}+b_{21}c_{21}x_{1}y_{1}y_{2}+a_{21}c_{22}x_{1}y_{1}y_{2}+\\
 & 2c_{21}c_{22}x_{3}y_{1}y_{2}+b_{21}c_{12}y{}_{2}^{2}-2b_{11}c_{22}y{}_{2}^{2}+b_{21}c_{22}x_{1}y{}_{2}^{2}+c_{22}^{2}x_{3}y{}_{2}^{2},\\
b_{12}= & -2b_{11}x_{1}+b_{21}x{}_{1}^{2}-c_{12}x_{3}+c_{22}x_{1}x_{3}+ty_{1},\\
a_{12}= & -2a_{11}x_{1}+a_{21}x{}_{1}^{2}-c_{11}x_{3}+c_{21}x_{1}x_{3}-ty_{2}.
\end{align*}

\vspace{3pt}

\noindent \underline{\bf $s$-chart}: Similarly to the $x_{1}$-chart,
we have only to describe $\mathfrak{U}_{\mA}^{14}\cap\{s=1\}$ as
follows: 
\begin{align*}
 & x=y^{\sharp},\\
 & t=1/3\,\beta(x,y).
\end{align*}

\vspace{3pt}

\noindent \underline{\bf $t$-chart}: Similarly to the $x_{1}$-chart,
we have only to describe $\mathfrak{U}_{\mA}^{14}\cap\{t=1\}$ as
follows: 
\begin{align*}
 & y=x^{\sharp},\\
 & s=1/3\,\beta(x,y).
\end{align*}

\vspace{3pt}
The above descriptions of charts show that the $x_{1}$-, $x_{2}$-,
$s$-, and $t$-charts of $\mathfrak{U}_{\mA}^{14}$ are isomorphic
to $\mA^{13}\times\mA^{1*}$. As for $s$- and $t$-charts, we also
refer to Corollary \ref{cor:affinespace}.

\vspace{3pt}

\noindent \underline{\bf $x_3$-chart}: Similarly to the $x_{1}$-chart,
it suffices to describe $\mathfrak{U}_{\mA}^{14}\cap\{x_{3}=1\}$.
When $x_{3}=1$, we may eliminate the coordinate $s$ by the 3rd entry
of the equation $sx=y^{\sharp}$ with $x_{3}=1$. Then we may verify
that $\mathfrak{U}_{\mA}^{14}\cap\{x_{3}=1\}$ is defined by the five
$4\times4$ Pfaffians of the following skew-symmetric matrix: {\small{}
\[
\left(\begin{array}{ccccc}
0 & t & a_{11}x_{1}+a_{12}x_{2}+c_{11} & b_{11}x_{1}+b_{12}x_{2}+c_{12} & x_{1}\\
 & 0 & a_{21}x_{1}-a_{11}x_{2}+c_{21} & b_{21}x_{1}-b_{11}x_{2}+c_{22} & x_{2}\\
 &  & 0 & y_{3} & -y_{2}\\
 &  &  & 0 & y_{1}\\
 &  &  &  & 0
\end{array}\right).
\]
}From this description, we see that {\small{}
\[
({\rm Sing}\,\mathfrak{U}_{\mA}^{14})\cap\{x_{3}\not=0\}=\{x_{3}\not=0,\,x_{1}=x_{2}=0,\,y=0,\,t=0,\,C=O\}.
\]
}{\small\par}

\vspace{3pt}

We define the following skew-symmetric matrix: {\small{}
\[
A_{y}=\left(\begin{array}{ccccc}
0 & s & a_{11}y_{1}+b_{11}y_{2}+y_{3} & a_{12}y_{1}+b_{12}y_{2} & c_{11}y_{1}+c_{12}y_{2}\\
 & 0 & a_{21}y_{1}+b_{21}y_{2} & -a_{11}y_{1}-b_{11}y_{2}+y_{3} & c_{21}y_{1}+c_{22}y_{2}\\
 &  & 0 & x_{3} & -x_{2}\\
 &  &  & 0 & x_{1}\\
 &  &  &  & 0
\end{array}\right).
\]
}{\small\par}

\vspace{3pt}

\noindent \underline{\bf $y_1$-chart}: Similarly to the $x_{3}$-chart,
We may verify that $\mathfrak{U}_{\mA}^{14}\cap\{y_{1}=1\}$ is defined
by the five $4\times4$ Pfaffians of the skew-symmetric matrix $A_{y}$
with $y_{1}=1$. From this description, we see that {\small{}
\begin{align*}
({\rm Sing}\,\mathfrak{U}_{\mA}^{14})\cap\{y_{1}\not=0\} & =\{y_{1}\not=0,\,x=0,\,s=0,\,M_{y}=O\}\\
 & =\left\{ y_{1}\not=0,\,x=0,\,s=y_{3}=0,\,\left(\begin{array}{cc}
a_{11} & b_{11}\\
a_{12} & b_{12}\\
c_{11} & c_{12}\\
c_{21} & c_{22}
\end{array}\right)\left(\begin{array}{c}
1\\
y_{2}
\end{array}\right)=\bm{o}\right\} .
\end{align*}
}{\small\par}

\vspace{3pt}

\noindent \underline{\bf $y_2$-chart}: Similarly to the $x_{3}$-chart,
We may verify that $\mathfrak{U}_{\mA}^{14}\cap\{y_{2}=1\}$ is defined
by the five $4\times4$ Pfaffians of the following skew-symmetric
matrix $A_{y}$ with $y_{2}=1$. From this description, we see that 

{\small{}
\begin{align*}
({\rm Sing}\,\mathfrak{U}_{\mA}^{14})\cap\{y_{2}\not=0\} & =\{y_{2}\not=0,\,x=0,\,s=0,\,M_{y}=O\}\\
 & =\left\{ y_{2}\not=0,\,x=0,\,s=y_{3}=0,\,\left(\begin{array}{cc}
a_{11} & b_{11}\\
a_{12} & b_{12}\\
a_{21} & b_{21}\\
c_{11} & c_{12}\\
c_{21} & c_{22}
\end{array}\right)\left(\begin{array}{c}
y_{1}\\
1
\end{array}\right)=\bm{o}\right\} .
\end{align*}
}{\small\par}

\vspace{3pt}

\noindent \underline{\bf $y_3$-chart}: We may easily verify that
the $y_{3}$-chart is contained in one of the above charts.

\vspace{10pt}

Let {\small{}
\[
\mathsf{S}_{x}:=\{x_{1}=x_{2}=0,\,y=0,\,s=0,\,t=0,\,C=O\}\simeq\mA^{7},
\]
}and $\mathsf{S}_{y}$ be the closure of the locus{\small{}
\[
\left\{ x=0,\,s=t=y_{3}=0,\,\left(\begin{array}{cc}
a_{11} & b_{11}\\
a_{12} & b_{12}\\
a_{21} & b_{21}\\
c_{11} & c_{12}\\
c_{21} & c_{22}
\end{array}\right)\left(\begin{array}{c}
y_{1}\\
y_{2}
\end{array}\right)=\bm{o}\right\} \setminus\{x=0,\,y=0,\,s=t=0\}.
\]
}{\small\par}

By a consideration with a little linear algebra, we see that {\small{}
\[
\mathsf{S}_{y}=\left\{ x=0,\,s=t=y_{3}=0,\,\rank\left(\begin{array}{cc}
a_{11} & b_{11}\\
a_{12} & b_{12}\\
a_{21} & b_{21}\\
c_{11} & c_{12}\\
c_{21} & c_{22}\\
-y_{2} & y_{1}
\end{array}\right)\leq1\right\} ,
\]
}which is the affine cone over $\mP^{1}\times\mP^{5}$ and hence is
$7$-dimensional.

Investigating the 9 equations of $\mathfrak{U}_{\mA}^{14}$, we see
that 
\begin{align*}
 & \mathsf{S}_{x}\cap\{x_{3}\not=0\}=({\rm Sing}\,\mathfrak{U}_{\mA}^{14})\cap\{x_{3}\not=0\},\\
 & \mathsf{S}_{y}\cap\{y_{i}\not=0\}=({\rm Sing}\,\mathfrak{U}_{\mA}^{14})\cap\{y_{i}\not=0\}\,(i=1,2),
\end{align*}

Thus, from the above descriptions of the charts, we can describe the
singularities of $\mathfrak{U}_{\mA}^{14}$ as follows:
\begin{prop}
\label{prop:The-singular-locus} The open subset $\mathfrak{U}_{\mA}^{14}\setminus\{x=y=0,s=t=0\}$
is $14$-dimensional and irreducible. Its singular locus is equal
to the $7$-dimensional locus
\[
\Delta:=(\mathsf{S}_{x}\cup\mathsf{S}_{y})\setminus\{x=y=0,s=t=0\},
\]
and it has $c({\rm G}(2,5))$-singularities along $\Delta$, where
we call a singularity isomorphic to the vertex of the cone over ${\rm G}(2,5)$
a \textup{$c({\rm G}(2,5))$-singularity }(it is easy to verify that
this is a terminal singularity by computing the discrepancy of the
exceptional divisor of the simple blowing up at the point).
\end{prop}

We will show that $\mathfrak{U}_{\mA}^{14}$ itself is $14$-dimensional
and irreducible in Proposition \ref{prop:Jerry14}, and determine
the singularities of $\mathfrak{U}_{\mA}^{14}$ along $\{x=0,y=0,s=t=0\}$
in Proposition \ref{prop:UTerm}.

\subsection{Gorensteinness and $9\times16$ graded minimal free resolution of
the ideal of $\mathfrak{U}_{\mathbb{A}}^{14}$\label{subsec:-is-Gorenstein}\protect 
}\begin{prop}
\label{prop:Jerry14} Let $S_{\mathfrak{U}}$ be the polynomial ring
over $\mC$ whose variables are $s,t$ and the entries of $A,B,C,x,y$.
Let $I_{\mathfrak{U}}$ be the ideal of the polynomial ring $S_{\mathfrak{U}}$
generated by the $9$ equations of $\mathfrak{U}_{\mathbb{A}}^{14}$.
Set $R_{\mathfrak{U}}:=S_{\mathfrak{U}}/I_{\mathfrak{U}}$. The following
assertions hold\,$:$ 

\begin{enumerate}[$(1)$]

\item We give nonnegative weights for the coordinates of $S_{\mathfrak{U}}$
such that all the equations of $\mathfrak{U}_{\mA}^{14}$ are weighted
homogeneous, and we denote by $w(*)$ the weight on the monomial $*$.
We denote by $\mP$ the corresponding weighted projective space, and
by $\mathfrak{U}_{\mP}\subset\mP$ the weighted projectivization of
$\mathfrak{U}_{\mA}^{14}$, where we allow some coordinates being
nonzero constants (thus $\dim\mathfrak{U}_{\mP}$ could be less than
$13$). We set
\begin{align*}
\delta & =2w(x_{1}x_{2})-w(x_{3})-w(y_{1}y_{2})+5w(y_{3}).
\end{align*}

\begin{enumerate}[$({1}\text{-}1)$]

\item It holds that 
\[
\omega_{\mP}=\sO_{\mP}(-4w(x_{1}x_{2})+5w(y_{1}y_{2})+4w(x_{3})-14w(y_{3})).
\]

\item The ideal $I_{\mathfrak{U}}$ has the following graded minimal
$S_{\mathfrak{U}}$-free resolution
\begin{equation}
0\leftarrow P_{0}\leftarrow P_{1}\leftarrow P_{2}\leftarrow P_{3}\leftarrow P_{4}\leftarrow0,\ \text{where}\label{eq:minfreeU}
\end{equation}
\begin{align*}
P_{0}= & S,\\
P_{1}= & S(-w(x_{1}y_{3}))\oplus S(-w(x_{2}y_{3}))\\
\oplus & S(-w(sx_{1}))\oplus S(-w(sx_{2}))\oplus S(-w(sx_{3}))\\
\oplus & S(-w(ty_{1}))\oplus S(-w(ty_{2}))\oplus S(-w(ty_{3}))\\
\oplus & S(-w(st)),\\
P_{2}= & S(-w(sx_{1}y_{3}))\oplus S(-w(sx_{2}y_{3}))\oplus S(-w(tx_{1}y_{3}))\oplus S(-w(tx_{2}y_{3}))\\
\oplus & S(-w(stx_{1}))\oplus S(-w(stx_{2}))\oplus S(-w(stx_{3}))\\
\oplus & S(-w(sty_{1}))\oplus S(-w(sty_{2}))\oplus S(-w(sty_{3}))\\
\oplus & S(-w(sx_{1}x_{2}))\oplus S(-w(sx_{1}x_{3}))\oplus S(-w(sx_{2}x_{3}))\\
\oplus & S(-w(ty_{1}y_{2}))\oplus S(-w(ty_{1}y_{3}))\oplus S(-w(ty_{2}y_{3}))\\
P_{3}= & S(-(\delta-w(x_{1}y_{3})))\oplus S(-(\delta-w(x_{2}y_{3})))\\
\oplus & S(-(\delta-w(sx_{1})))\oplus S(-(\delta-w(sx_{2})))\oplus S(-(\delta-w(sx_{3})))\\
\oplus & S(-(\delta-w(ty_{1})))\oplus S(-(\delta-w(ty_{2})))\oplus S(-(\delta-w(ty_{3})))\\
\oplus & S(-(\delta-w(st))),\\
P_{4}= & S(-\delta).
\end{align*}

\item It holds that 
\begin{equation}
\omega_{\mathfrak{U}_{\mathbb{P}}}=\sO_{\mathfrak{U}_{\mathbb{P}}}(-2w(x_{1}x_{2})+4w(y_{1}y_{2})+3w(x_{3})-9w(y_{3})).\label{eq:omegaP}
\end{equation}

\end{enumerate}

\item $I_{\mathfrak{U}}$ is a Gorenstein ideal of codimension $4$. 

\item $\mathfrak{U}_{\mathbb{A}}^{14}$ is irreducible and reduced,
thus $I_{\mathfrak{U}}$ is a prime ideal.

\item $\mathfrak{U}_{\mathbb{A}}^{14}$ is normal. 

\end{enumerate}
\end{prop}

\begin{proof}
We may compute the $S_{\mathfrak{U}}$-free resolution (\ref{eq:minfreeU})
of $I_{\mathfrak{U}}$ by \textsc{Singular} \cite{DGPS}. For the
remaining assertions, the proof of \cite[Prop.4.8]{Tak4} works verbatim. 
\end{proof}

\subsection{Factoriality of $\mathfrak{U}_{\mathbb{A}}^{14}$\label{subsec:FactorialityU}}

We denote by $\bar{x_{1}}$ the image of $x_{1}$ in $R_{\mathfrak{U}}$. 
\begin{lem}
\label{lem:primeU} $\mathfrak{U}_{\mathbb{A}}^{14}\cap\{x_{1}=0\}$
is irreducible and reduced, and is normal. In particular, the element
$\bar{x_{1}}\in R_{\mathfrak{U}}$ is a prime element.
\end{lem}

\begin{proof}
By \cite[Lem.3.6]{Tak5}, it suffices to show that ${\rm Sing}\,(\mathfrak{U}_{\mathbb{A}}^{14}\cap\{x_{1}=0\})$
has codimension $\geq2$ in $\mathfrak{U}_{\mathbb{A}}^{14}\cap\{x_{1}=0\}$.
By the descriptions of charts as in Subsection \ref{subsec:Descriptrion-of-charts},
we see that $\mathfrak{U}_{\mathbb{A}}^{14}\cap\{x_{1}=0\}$ is smooth
on the $x_{2}$- and the $t$-charts, and is isomorphic to a hypersurface
\begin{align*}
 & \{-a_{11}c_{11}y_{1}^{2}-a_{12}c_{21}y{}_{1}^{2}-b_{11}c_{11}y_{1}y_{2}-a_{11}c_{12}y_{1}y_{2}-b_{12}c_{21}y_{1}y_{2}\\
 & -a_{12}c_{22}y_{1}y_{2}-b_{11}c_{12}y{}_{2}^{2}-b_{12}c_{22}y{}_{2}^{2}+c_{11}y_{1}y_{3}+c_{12}y_{2}y_{3}=0\}
\end{align*}
on the $s$-chart, whose singular locus has codimension $\geq2$ in
$\mathfrak{U}_{\mathbb{A}}^{14}\cap\{x_{1}=0\}$. The locus $\mathfrak{U}_{\mathbb{A}}^{14}\cap\{x_{1}=x_{2}=0,s=t=0\}$
is the union of the following two loci (i) and (ii):

\begin{enumerate}[(i)]

\item $\{x=0,s=t=0,y^{\sharp}=0\}$, which is a fibration of relative
dimension $7$ over the affine cone of $\mP^{1}\times\mP^{2}$, hence
is $11$-dimensional.

\item $\{\rank\left(\begin{array}{ccc}
c_{11} & c_{21} & -y_{2}\\
c_{12} & c_{22} & y_{1}
\end{array}\right)\leq1,F_{3}=0\},$ which is easily to be seen of codimension $2$ in $\mathfrak{U}_{\mathbb{A}}^{14}\cap\{x_{1}=0\}$
($F_{3}$ is defined in the proof of (1)).

\end{enumerate}

Therefore we have shown that ${\rm Sing}\,(\mathfrak{U}_{\mathbb{A}}^{14}\cap\{x_{1}=0\})$
has codimension $\geq2$ in $\mathfrak{U}_{\mathbb{A}}^{14}\cap\{x_{1}=0\}$.
\end{proof}
\begin{prop}
\label{prop:UFD-U}The affine coordinate ring $R_{\mathfrak{U}}$
of $\mathfrak{U}_{\mathbb{A}}^{14}$ is a UFD. 
\end{prop}

\begin{proof}
Using the description of the $x_{1}$-chart as in Subsection \ref{subsec:Descriptrion-of-charts}
and Lemma \ref{lem:primeU}, the proof of \cite[Prop.4.9]{Tak4} works
verbatim.
\end{proof}
The following corollary can be proved in the same way as the proof
of \cite[Cor.4.10]{Tak4}:
\begin{cor}
\label{UFactorial} Let $\mathfrak{U}_{\mP}^{13}$ be the weighted
projectivization of $\mathfrak{U}_{\mA}^{14}$ with some positive
weights on the coordinates. The following assertions hold:

\begin{enumerate}

\item Any prime Weil divisor on $\mathfrak{U}_{\mP}^{13}$ is the
intersection between $\mathfrak{U}_{\mP}^{13}$ and a weighted hypersurface.
In particular, $\mathfrak{U}_{\mP}^{13}$ is $\mQ$-factorial and
has Picard number one. 

\item Let $X$ be a quasi-smooth threefold such that $X$ is a codimension
$10$ weighted complete intersection in $\mathfrak{U}_{\mP}^{13}$,
i.e., there exist ten weighted homogeneous polynomials $G_{1},\dots,G_{10}$
such that $X=\mathfrak{U}_{\mP}^{13}\cap\{G_{1}=0\}\cap\dots\cap\{G_{10}=0\}$.
Assume moreover that $\{x_{1}=0\}\cap X$ is a prime divisor. Then
any prime Weil divisor on $X$ is the intersection between $X$ and
a weighted hypersurface. In particular, $X$ is $\mQ$-factorial and
has Picard number one.\end{enumerate}
\end{cor}

\subsection{$\mP^{1}\times\mP^{1}\times\mP^{1}$-fibration \label{subsec:-fibration}}

We denote by $\mA_{\mathsf{\mathtt{P}}}\simeq\mA^{10}$ the affine
space whose coordinates are the entries of the matrices $A,B,C$.
Note that any equation of $\mathfrak{U}_{\mA}^{14}$ is of degree
two if we regard the entries of $A,B,C$ as constants. Therefore,
considering the variables of the equations of $\mathfrak{U}_{\mA}^{14}$
except the entries of $A,B,C$ as projective coordinates, we obtain
a $13$-dimensional quasi-projective variety with the same equation
as $\mathfrak{U}_{\mA}^{14}$. We denote this $13$-dimensional variety
by $\widehat{\mathfrak{U}}$, and by $\rho_{\mathfrak{U}}\colon\widehat{\mathfrak{U}}\to\mA_{\mathtt{P}}$
the natural projection. Note that the ${\rm SL_{2}}\times{\rm SL}_{2}$-action
on $\mathfrak{U}_{\mA}^{14}$ defined as in Subsection \ref{subsec:SL2SL2}
induces those on $\widehat{\mathfrak{U}}$ and those on $\mA_{\mathtt{P}}$,
and the natural projection $\rho_{\mathfrak{U}}\colon\widehat{\mathfrak{U}}\to\mA_{\mathtt{P}}$
is equivariant with respect to these actions. Using these group actions,
detail analysis of $\rho_{\mathfrak{U}}\colon\widehat{\mathfrak{U}}\to\mA_{\mathtt{P}}$
is possible. We, however, omit the details since they are lengthy.
Instead, we will describe the restriction of $\rho_{\mathfrak{U}}$
to the subvariety \textbf{$\mathfrak{S}_{\mA}^{8}$ }defined in Section
\ref{sec:KeyEquation8}. Here we only state the following, which can
be shown immediately by Corollary \ref{cor:P1P1P1}:
\begin{prop}
Let $\Delta_{\rho_{\mathfrak{U}}}$ be the closed subset of $\mA_{\mathtt{P}}$
defined by $D_{\beta}$ as in Remark \ref{rem:Dbeta}. The $\rho_{\mathfrak{U}}$-fiber
over a point outside $\Delta_{\rho_{\mathfrak{U}}}$ is isomorphic
to $\mP\times\mP^{1}\times\mP^{1}.$ 
\end{prop}

\subsection{Terminal singularities\label{subsec:Terminal-singularities}}

By Proposition \ref{prop:The-singular-locus} and the result in Subsection
\ref{subsec:-fibration}, we can show the following in the same way
as \cite[Prop.4.11]{Tak4}:
\begin{prop}
\label{prop:UTerm} The variety $\mathfrak{U}_{\mA}^{14}$ has only
terminal singularities with the following descriptions: 
\begin{enumerate}
\item The singularities along the $7$-dimensional locus $(\mathsf{S}_{x}\cup\mathsf{S}_{y})\setminus\{x=0,y=0,s=t=0\}$
are $c({\rm G}(2,5))$-singularities (cf.~Proposition \ref{prop:The-singular-locus}). 
\item There exists a primitive $K$-negative divisorial extraction $f\colon\widetilde{\mathfrak{U}}\to\mathfrak{U}_{\mA}^{14}$
such that 
\begin{enumerate}
\item singularities of $\widetilde{\mathfrak{U}}$ are only $c({\rm G}(2,5))$-singularities
along the strict transforms of $\mathsf{S}_{x}\cup\mathsf{S}_{y}$,
and 
\item for the $f$-exceptional divisor $E_{\mathfrak{U}}$, the morphism
$f|_{E_{\mathfrak{U}}}$ can be identified with $\rho_{\mathfrak{U}}\colon\widehat{\mathfrak{U}}^{13}\to\mA_{{\rm B}}^{10}$
as in Subsection \ref{subsec:-fibration}. 
\end{enumerate}
\end{enumerate}
\end{prop}

\subsection{An example of a $\mQ$-Fano threefold\label{subsec:An-example-of fano 1}}

We consider the weighted projectivization of $\mathfrak{U}_{\mathbb{A}}^{14}$
\begin{align*}
 & \mathfrak{U}_{\mP}^{13}\subset\mP(1^{15},2^{2},3)
\end{align*}
by the weights on the coordinates as in Example \ref{exa:In-Subsection-,WT U}.
In this subsection, we show the following:
\begin{thm}
\label{thm:20544Ex}Let $L_{1},\dots,L_{10}$ be general forms of
weight $1$ in $\mP(1^{15},2^{2},3)$. The subscheme
\[
X:=\mathfrak{U}_{\mP}^{13}\cap L_{1}\cap\cdots\cap L_{10}
\]
is a prime $\mQ$-Fano $3$-fold of No.20544 and of anticanonical
codimension $4$. 
\end{thm}

\begin{proof}
In this proof, we denote $\mP(1^{15},2^{2},3)$ by $\mP$ for simplicity.
Using the equations of $\mathfrak{U}_{\mP}^{13}$, we may easily verify
that ${\rm Bs}|\sO_{\mP}(1)|\cap\mathfrak{U}_{\mP}^{13}$ consists
of the $s$-point and the $t$-point, where the $*$-point for a coordinate
$*$ of $\mP$ means the point of $\mP$ such that all the coordinates
except $*$ are zero. Note that $\dim{\rm Sing\,}\mathfrak{U}_{\mathbb{A}}^{14}$
is less than the codimension of $X$ in $\mathfrak{U}_{\mathbb{P}}^{13}$.
Therefore, by the Bertini theorem, we see that $X$ is a smooth $3$-fold
outside the $s$-point and the $t$-point. Computing the linear parts
of the equations of $X$ at each of the $s$-point and the $t$-point
(cf. LPC in \cite[Subsec.2.8.1]{Tak7}), we see that $X$ has a $1/3(1,1,2)$-singularity
at the $s$-point and a $1/2(1,1,1)$-singularity at the $t$-point.
By Proposition \ref{prop:Jerry14} (1-2), we have the following $S_{\mathfrak{U}}$-free
minimal resolution of $R_{\mathfrak{U}}$ which is graded with respect
to the weights on the variables given in Example \ref{exa:In-Subsection-,WT U}:
\begin{align*}
 & 0\leftarrow R_{\mathfrak{U}}\leftarrow S_{\mathfrak{U}}\leftarrow S_{\mathfrak{U}}(-3)^{\oplus4}\oplus S_{\mathfrak{U}}(-4)^{\oplus4}\oplus S_{\mathfrak{U}}(-5)\\
 & \leftarrow S_{\mathfrak{U}}(-4)\oplus S_{\mathfrak{U}}(-5)^{\oplus7}\oplus S_{\mathfrak{U}}(-6)^{\oplus7}\oplus S_{\mathfrak{U}}(-7)\\
 & \leftarrow S_{\mathfrak{U}}(-6)\oplus S_{\mathfrak{U}}(-7)^{\oplus4}\oplus S_{\mathfrak{U}}(-8)^{\oplus4}\leftarrow S_{\mathfrak{U}}(-11)\leftarrow0.
\end{align*}
From this, we see that $-K_{X}=\sO_{X}(1)$ and $(-K_{X})^{3}=31/6$.
It remains to show that $\rho(X)=1$. By Corollary \ref{UFactorial},
it suffices to check that $X\cap\{x_{1}=0\}$ is a prime divisor.
This follows by the Bertini theorem since $\mathfrak{U}_{\mP}^{13}\cap\{x_{1}=0\}$
is normal by Lemma \ref{lem:primeU}, and the fact that ${\rm Bs}|\sO_{\mP}(1)|\cap\mathfrak{U}_{\mP}^{13}$
consists of the $s$-point and the $t$-point. Hence $X$ is a prime
$\mQ$-Fano 3-fold of No.20544, and by \cite[Subsec.2.5, Proof  of  Thm.1  (1)]{Tak7},
$X$ is of anti-canonical codimension $4$. 
\end{proof}
\begin{rem}
In a future work, we will describe via $\mathfrak{U}_{\mP}^{13}$
the Sarkisov link for a prime $\mQ$-Fano threefold of No.20544 starting
from the weighted blow-up at the unique $\text{1/3(1,1,2)}$-singularity. 
\end{rem}

\section{\textbf{A specialization of $\mathfrak{U}_{\mathbb{A}}^{14}$--an
affine variety $\mathfrak{S}_{\mA}^{8}$ with an ${\rm SL}_{2}$-action}--\label{sec:KeyEquation8}}

%Notation change:Note 39 p47 

In this section, we define a closed subvariety $\mathfrak{S}_{\mathbb{A}}^{8}$
of the affine variety $\mathfrak{U}_{\mathbb{A}}^{14}$ and investigate
its properties. The whole story in this section is very similar to
that of Section \ref{sec:KeyEquation14}, hence we will not write
down it fully. 

As we will see in Subsection \ref{subsec:The-case-A3A4}, $\mathfrak{S}_{\mathbb{A}}^{8}$
is isomorphic to a subvariety of the $G_{2}^{(4)}$-cluster variety.
We will discuss there a prime $\mQ$-Fano threefold of anti-canonical
codimension 4 which is obtained as a weighted complete intersection
of a weighted projectivization of $\mathfrak{S}_{\mathbb{A}}^{8}$.

\subsection{Definition}
\begin{defn}
We set 
\[
\mathfrak{S}_{\mA}^{8}:=\mathfrak{U}_{\mathbb{A}}^{14}\cap\left\{ b_{11}=a_{12},b_{21}=-a_{11},\left(\begin{array}{cc}
c_{11} & c_{12}\\
c_{21} & c_{22}
\end{array}\right)=\left(\begin{array}{cc}
1 & 0\\
0 & 1
\end{array}\right)\right\} .
\]
\end{defn}

\begin{rem}
We recall that the matrices $M_{x}$ and $M_{y}$ are defined in Subsection
\ref{subsec:DefinitionU}. We note that the condition on $\mathfrak{S}_{\mA}^{8}$
is equivalent to that $M_{x}$ coincides with $M_{y}$ by replacing
$x_{1},x_{2},x_{3}$ with $y_{1},y_{2},y_{3}$ respectively.
\end{rem}

Now we introduce a different notation for coordinates of $\mathfrak{S}_{\mA}^{8}$
as follows, with which we derive another presentation of the equations
of $\mathfrak{S}_{\mA}^{8}$ (Proposition \ref{prop:NewS8}) such
that a group action on $\mathfrak{S}_{\mathbb{A}}^{8}$ will be more
visible (Subsection \ref{subsec:-action}):
\begin{align*}
 & P=\left(\begin{array}{ccc}
d_{2} & d_{1} & 1\\
d_{1} & d_{0} & 0\\
1 & 0 & 0
\end{array}\right),Q=\left(\begin{array}{ccc}
-d_{3} & -d_{2} & 0\\
-d_{2} & -d_{1} & 1\\
0 & 1 & 0
\end{array}\right),\\
 & U:=\left(\begin{array}{cc}
u_{1} & u_{2}\\
u_{3} & u_{4}
\end{array}\right):=\left(\begin{array}{cc}
y_{1} & -x_{1}\\
y_{2} & -x_{2}
\end{array}\right),\\
 & V:=\left(\begin{array}{ccc}
v_{2} & v_{1} & v_{0}\\
-v_{3} & -v_{2} & -v_{1}
\end{array}\right):=\left(\begin{array}{ccc}
y_{3} & x_{3} & t\\
-s & -y_{3} & -x_{3}
\end{array}\right).
\end{align*}
We also set
\begin{align*}
 & \widehat{U}:=\left(\begin{array}{ccc}
u_{1}^{2} & u_{1}u_{2} & u_{2}^{2}\\
2u_{1}u_{3} & u_{1}u_{4}+u_{2}u_{3} & 2u_{2}u_{4}\\
u_{3}^{2} & u_{3}u_{4} & u_{4}^{2}
\end{array}\right),\,\widehat{U}^{\dagger}:=\left(\begin{array}{ccc}
u_{4}^{2} & -u_{2}u_{4} & u_{2}^{2}\\
-2u_{3}u_{4} & u_{1}u_{4}+u_{2}u_{3} & -2u_{1}u_{2}\\
u_{3}^{2} & -u_{1}u_{3} & u_{1}^{2}
\end{array}\right),\\
 & D:=\left(\begin{array}{ccc}
d_{2} & d_{1} & d_{0}\\
-d_{3} & -d_{2} & -d_{1}
\end{array}\right).
\end{align*}
By a straightforward calculation, we can check the following:
\begin{prop}
\label{prop:NewS8}The affine variety $\mathfrak{S}_{\mA}^{8}$ is
defined by the following equations:
\end{prop}

\begin{equation}
UV=D\widehat{U},\,\wedge^{2}V=(\wedge^{2}D){\empty}^{t}\!\widehat{U}^{\dagger}.\label{eq:S8}
\end{equation}

\begin{rem}
If $\det U\not=0$, then $\wedge^{2}V=(\wedge^{2}D){\empty}^{t}\!\widehat{U}^{\dagger}$
can be derived from $UV=D\widehat{U}$ by the Cauchy-Binet formula.
\end{rem}

\subsection{$(({\rm SL}_{2}\rtimes{\rm SL_{2})}\times(\mC^{*})^{2})$-action\label{subsec:-action}}

We denote by $S_{\mathfrak{S}}$ the polynomial ring over the field
$\mC$ with the entries of $U,V,D$ as the variables. We consider
the affine space $\mathbb{A}_{\mathfrak{S}}^{12}$ with the coordinate
ring $S_{\mathfrak{S}}$. We denote by $\mA_{D}^{4}$ the affine space
with the entries of $D$ as the coordinates.

From the presentation of the equation of $\mathfrak{S}_{\mA}^{8}$
as in Proposition \ref{prop:NewS8}, we can immediately read off an
$\left(({\rm SL_{2}}\rtimes{\rm SL_{2}})\times(\mC^{*})^{2}\right)$-action
on $\mathfrak{S}_{\mA}^{8}$ as follows: 
\begin{prop}
\label{prop:S8Orbit}The following assertions hold:

\begin{enumerate}

\item 

\begin{enumerate}[({1}-1)]

\item The variety $\mathfrak{S}_{\mA}^{8}$ is preserved by the following
actions of the two groups ${\rm SL_{2}^{I}}$, ${\rm SL_{2}^{II}}$
isomorphic to ${\rm SL_{2}}$, and they define an action of the group
${\rm SL_{2}^{II}}\rtimes{\rm SL_{2}^{I}}$ on $\mathfrak{S}_{\mathbb{A}}^{8}$:

\begin{align*}
 & \text{For}\ g\in{\rm SL_{2}^{I}},U\mapsto gUg^{-1},\,V\mapsto gV\widehat{g}^{-1},\,D\mapsto gD\widehat{g}^{-1},\\
 & \text{For}\ h\in{\rm SL_{2}^{II}},U\mapsto Uh,\,V\mapsto h^{-1}V\widehat{h},\,D\mapsto D,
\end{align*}
where the definitions of $\widehat{g}$ and $\widehat{h}$ for $g$
and $h$ respectively are similar to that of $\widehat{U}$ for $U$.

\item The affine space $\mathbb{A}_{\mathfrak{S}}^{12}$ has the
$(\mC^{*})^{2}$-action preserving $\mathfrak{S}_{\mA}^{8}$ defined
by

\[
U\mapsto\alpha U,\,D\mapsto\beta D,\,V\mapsto\alpha\beta V,
\]
where $\alpha,\beta\in\mC^{*}$. 

\end{enumerate}

These induce an $\left(({\rm SL_{2}^{II}}\rtimes{\rm SL_{2}^{I}})\times(\mC^{*})^{2}\right)$-action
on $\mathfrak{S}_{\mA}^{8}$. 

\item The induced $({\rm SL_{2}^{I}}\times(\mC^{*})^{2})$-action
on $\mA_{D}^{4}$ has the following orbits:

\begin{enumerate}[(a)]

\item $\{0\}$.

\item The complement of $\{0\}$ in the cone over the twisted cubic
$\gamma_{D}$ defined by $\wedge^{2}D=\bm{0}$.

\item The complement of the orbits as described in $(a)$ and $(b)$
in the cone over the tangential scroll of the twisted cubic $\gamma_{D}$,
where the equation of the tangential scroll is 
\begin{equation}
3d_{1}^{2}d_{2}^{2}-4d_{1}^{3}d_{3}-4d_{0}d_{2}^{3}+6d_{0}d_{1}d_{2}d_{3}-d_{0}^{2}d_{3}^{2}=0.\label{eq:tangential}
\end{equation}

\item The complement of the orbits as described in $(a)$, $(b)$
and $(c)$ in \textup{$\mA_{D}^{4}$.}

\end{enumerate}

\end{enumerate}
\end{prop}

\begin{proof}
For (1), we only note that $\widehat{U}$ is mapped to $\widehat{g}\widehat{U}\widehat{g}^{-1}$
by the ${\rm SL_{2}^{I}}$-action. The assertion $(2)$ is well-known
to be true.
\end{proof}
\begin{rem}
By an explicit calculation, we may easily verify that the rank of
the matrix associated to the bi-linear trace $\beta$ for $\mathfrak{S}_{\mA}^{8}$
is $0,1,2,$ or $3$ if and only if the Proposition \ref{prop:S8Orbit}
(2) (a), (b), (c), or (d) holds respectively.
\end{rem}

\subsection{$\mP^{1}\times\mP^{1}\times\mP^{1}$-fibration\label{subsec:-fibration-8dim}}

We denote by $\mP_{\mathtt{F}}^{7}$ the projective space whose coordinates
are the entries of the matrices $U,V$. Note that any equation of
$\mathfrak{S}_{\mA}^{8}$ is of degree two if we regard the entries
of $D$ as constants. Therefore, considering the variables of the
equations of $\mathfrak{S}_{\mA}^{8}$ except the entries of $D$
as projective coordinates, we obtain a $7$-dimensional variety in
$\mA_{{\rm D}}^{4}\times\mP_{\mathtt{F}}^{7}$. We denote this $7$-dimensional
variety by $\widehat{\mathfrak{S}}^{7}$, and by $\rho_{\mathfrak{S}}\colon\widehat{\mathfrak{S}}^{7}\to\mA_{D}^{4}$
the natural projection. Note that the $({\rm SL}_{2}\times(\mC^{*})^{2})$-action
on $\mathfrak{S}_{\mA}^{8}$ defined as in Subsection \ref{subsec:-action}
induces those on $\widehat{\mathfrak{S}}^{7}$ and those on $\mA_{{\rm D}}$,
and the natural projection $\rho_{\mathfrak{S}}\colon\widehat{\mathfrak{S}}^{7}\to\mA_{{\rm D}}^{4}$
is equivariant with respect to these actions. Using these group actions,
we will describe $\rho_{\mathfrak{S}}\colon\widehat{\mathfrak{S}}^{7}\to\mA_{D}^{4}$.
\begin{lem}
The following two assertions hold:

\begin{enumerate}

\item In the projective $7$-space with coordinates $x_{ij}\,(i=0,1,\,0\leq j\leq3)$,
let 
\[
\mathsf{S}:=\left\{ \rank\left(\begin{array}{cccc}
x_{00} & x_{01} & x_{02} & x_{03}\\
x_{10} & x_{11} & x_{12} & x_{13}
\end{array}\right)\leq1\right\} ,
\]
which is nothing but the image of the Segre embedding of $\mP^{1}\times\mP^{3}$.
Let $\mathsf{Q}\subset\mP^{3}$ be a cone over a smooth conic. The
$3$-fold $\mP^{1}\times\mathsf{Q}$ which is embedded in $\mP^{7}$
by the restriction of the Segre embedding of $\mP^{1}\times\mP^{3}$
is projectively equivalent to {\footnotesize{}
\begin{equation}
\left\{ \rank\left(\begin{array}{cccc}
x_{00} & x_{01} & x_{02} & x_{03}\\
x_{10} & x_{11} & x_{12} & x_{13}
\end{array}\right)\leq1,\,x_{03}^{2}=x_{01}x_{02},\,x_{13}^{2}=x_{11}x_{12},\,x_{03}x_{13}=x_{02}x_{11}\right\} .\label{eq:P1Q}
\end{equation}
}The projective variety defined by the equation (\ref{eq:P1Q}) is
a sextic del Pezzo $3$-fold with $A_{1}$ singularities along the
$(x_{00}:x_{10})$-line.

\item In the projective $7$-space with coordinates $u_{i}\,(1\leq i\leq4),v_{j}\,(0\leq j\leq3)$,
the variety defined by {\small{}
\begin{equation}
\left\{ \mathrm{rk}\left(\begin{array}{cccc}
u_{3} & v_{3} & v_{2} & v_{1}\\
u_{4} & v_{2} & v_{1} & v_{0}
\end{array}\right)\leq1,u_{1}v_{2}-u_{2}v_{3}=u_{3}^{2},u_{1}v_{1}-u_{2}v_{2}=u_{3}u_{4},u_{1}v_{0}-u_{2}v_{1}=u_{4}^{2}\right\} \label{eq:P111}
\end{equation}
}is a sextic del Pezzo $3$-fold with $A_{2}$ singularities along
the $(u_{1}:u_{2})$-line. This variety is isomorphic to the projective
$3$-fold $\mP^{1,1,1}$ defined in \cite[Def.6.6]{Fuk}.

\end{enumerate}
\end{lem}

\begin{proof}
(1). Let $s\colon\mP^{1}\times\mP^{3}\to\mathsf{S}$ be the Segre
embedding defined by $x_{ij}=p_{i}q_{j}$, where $p_{i}\,(i=0,1)$
and $q_{j}\,(0\leq j\leq3)$ are the coordinates of $\mP^{1}$ and
$\mP^{3}$ respectively. Then the inverse image of (\ref{eq:P1Q})
is $\{q_{3}^{2}=q_{1}q_{2}\}\subset\mP^{1}\times\mP^{3}$, which is
the product of $\mP^{1}$ and the cone over a smooth conic in $\mP^{3}$.
Therefore we have shown the former assertion. 

Now we show the latter assertion. By the equation (\ref{eq:P1Q}),
we see that the projective variety defined by the equation (\ref{eq:P1Q})
has $A_{1}$ singularities along the $(x_{00}:x_{10})$-line. We can
compute the minimal free resolution of the structure sheaf of the
variety defined by (\ref{eq:P1Q}) as follows:
\[
\sO(-2)^{\oplus9}\leftarrow\sO(-3)^{\oplus16}\leftarrow\sO(-4)^{\oplus9}\leftarrow\sO(-6).
\]
From this, we see that the projective variety defined by the equation
(\ref{eq:P1Q}) is a sextic del Pezzo $3$-fold.

(2). The former assertion can be shown similarly to the proof of the
latter assertion of (1). The latter assertion follows since $\mP^{1,1,1}$
is also a sextic del Pezzo $3$-fold with $A_{2}$ singularities along
a line by \cite[Rem.6.7]{Fuk} and there is only one such a sextic
del Pezzo $3$-fold (see \cite[(si31i)]{Fuj}).
\end{proof}
\begin{prop}
\label{prop:P1P1P1FibS8}Let $\mathsf{p}$ be a point of $\mA_{D}^{4}$
and $F_{\mathsf{p}}$ the $\rho_{\mathfrak{S}}$-fiber over $\mathsf{p}.$
We use the descriptions of the $({\rm SL_{2}^{I}}\times(\mC^{*})^{2})$-action
on $\mA_{D}^{4}$ as in Proposition \ref{prop:S8Orbit} $(a)$--$(d)$.

If $\mathsf{p}=0$, then $F_{\mathsf{p}}=\{UV=\bm{0},\wedge^{2}V=\bm{0}\}.$

If $\mathsf{p}$ belongs to the orbit as in $(b)$, then $F_{\mathsf{p}}\simeq\mP^{1,1,1}.$

If $\mathsf{p}$ belongs to the orbit as in $(c)$, then $F_{\mathsf{p}}\simeq\mP^{1}\times\mathsf{Q}$.

If $\mathsf{p}$ belongs to the orbit as in $(d)$, then $F_{\mathsf{p}}\simeq\mP^{1}\times\mP^{1}\times\mP^{1}.$
\end{prop}

\begin{proof}
If $\mathsf{p}=0$, then the description of $F_{\mathsf{p}}$ follows
from the equation of $\mathfrak{S}_{\mA}^{8}$.

By the $({\rm SL_{2}^{I}}\times(\mC^{*})^{2})$-action, we may choose
the point $\mathsf{p}$ as a special point as follows according to
the orbit to which $\mathsf{p}$ belongs:

If $\mathsf{p}$ belongs to the orbit as in $(b)$, we may assume
$\mathsf{p}$ is the $d_{0}$-point. Then the equation of $F_{\mathsf{p}}$
is {\small{}
\[
\left\{ \mathrm{rk}\left(\begin{array}{cccc}
u_{3} & v_{3} & v_{2} & v_{1}\\
u_{4} & v_{2} & v_{1} & v_{0}
\end{array}\right)\leq1,u_{1}v_{2}-u_{2}v_{3}=u_{3}^{2},u_{1}v_{1}-u_{2}v_{2}=u_{3}u_{4},u_{1}v_{0}-u_{2}v_{1}=u_{4}^{2}\right\} ,
\]
}which is nothing but $\mP^{1,1,1}.$

If $\mathsf{p}$ belongs to the orbit as in $(c)$, we may assume
$\mathsf{p}$ is the $d_{1}$-point. Then the equation of $F_{\mathsf{p}}$
coincides with (\ref{eq:P1Q}) by setting

\[
\left(\begin{array}{cccc}
x_{00} & x_{01} & x_{02} & x_{03}\\
x_{10} & x_{11} & x_{12} & x_{13}
\end{array}\right)=\left(\begin{array}{cccc}
u_{1} & v_{1}+2u_{4} & v_{3} & v_{2}+u_{3}\\
u_{2} & v_{0} & v_{2}-2u_{3} & v_{1}-u_{4}
\end{array}\right),
\]
hence $F_{\mathsf{p}}$ is isomorphic to $\mP^{1}\times\mathsf{Q}$.

Now assume that $\mathsf{p}$ belongs to the orbit as in $(d)$. Note
that the cone over the tangential scroll of the twisted cubic $\gamma_{D}$
as in (\ref{eq:tangential}) is the restriction of $D_{\beta}=0$
as in Remark \ref{rem:Dbeta}, which follows by a straightforward
calculation. Therefore, by Corollary \ref{cor:P1P1P1}, we see that
$F_{\mathsf{p}}\simeq\mP^{1}\times\mP^{1}\times\mP^{1}.$ 
\end{proof}

\subsection{Summary of properties of $\mathfrak{S}_{\mA}^{8}$\label{subsec:Summary-of-properties}}

In this subsection, we sum up the properties of $\mathfrak{S}_{\mA}^{8}$
corresponding to those of $\mathfrak{U}_{\mA}^{14}$ as in Subsections
\ref{subsec:-is-Gorenstein}, \ref{subsec:FactorialityU} and \ref{subsec:Terminal-singularities}
as follows:
\begin{prop}
\label{prop:summary8dim} The following assertions hold:

\begin{enumerate}

\item Let $S_{\mathfrak{S}}$ be the polynomial ring over $\mC$
whose variables are the entries of $U,V,D$. Let $I_{\mathfrak{S}}$
be the ideal of the polynomial ring $S_{\mathfrak{S}}$ generated
by the $9$ equations of $\mathfrak{S}_{\mathbb{A}}^{8}$. Set $R_{\mathfrak{S}}:=S_{\mathfrak{S}}/I_{\mathfrak{S}}$. 

\begin{enumerate}[({1}-1)]

\item $I_{\mathfrak{S}}$ is a prime, and Gorenstein ideal of codimension
four. 

\item $\mathfrak{S}_{\mathbb{A}}^{8}$ and $\mathfrak{S}_{\mathbb{A}}^{8}\cap\{u_{2}=0\}$
are normal. 

\item $R_{\mathfrak{S}}$ is a UFD. 
\end{enumerate}

\item The weighted projectivization $\mathfrak{S}_{\mP}^{7}$ of
$\mathfrak{S}_{\mA}^{8}$ with some positive weights on the coordinates
is $\mQ$-factorial and has Picard number one (we refer to Example
\ref{exa:We-consider-theEx20652} for an example of such a set of
weights). 

\item

There exists a primitive $K$-negative divisorial extraction $f\colon\widetilde{\mathfrak{S}}\to\mathfrak{S}_{\mA}^{8}$
such that $\widetilde{\mathfrak{S}}$ is smooth, and, for the $f$-exceptional
divisor $E_{\mathfrak{S}}$, the morphism $f|_{E_{\mathfrak{S}}}$
can be identified with $\rho_{\mathfrak{S}}\colon\widehat{\mathfrak{S}}^{7}\to\mA_{{\rm D}}^{4}$
as in Subsection \ref{subsec:-fibration-8dim}. In particular, the
variety $\mathfrak{S}_{\mA}^{8}$ has only terminal singularities. 

\end{enumerate}
\end{prop}

We omit the proof since we can prove this in the same (and simpler)
way as for the variety $\mathfrak{U}_{\mA}^{14}$.

\vspace{5pt}

As for examples of prime $\mQ$-Fano 3-folds obtained from $\mathfrak{S}_{\mA}^{8}$,
we see in Subsection \ref{subsec:The-case-A3A4} that they are actually
obtained from the $G_{2}^{(4)}$-variety. 

\section{\textbf{Another specialization of $\mathfrak{F}_{\mathbb{A}}^{22}$--an
affine variety $\mathfrak{Z}_{\mA}^{12}$ with an ${\rm SL_{3}}$-action--\label{sec:Z12}}}

In this section, we consider one more specialization of $\mathfrak{F}_{\mathbb{A}}^{22}$
as follows:
\begin{defn}
\label{def:We-set-Z12}We set 
\[
\mathfrak{Z}_{\mA}^{12}:=\mathfrak{F}_{\mathbb{A}}^{22}\cap\left\{ P=E,\,{\rm Tr}\,Q=0\right\} ,
\]
 where $E$ is the $3\times3$ identity matrix, and ${\rm Tr}\,Q$
is the trace of the matrix $Q=(q_{ij})$.
\end{defn}

Actually, the equations of $\mathfrak{Z}_{\mA}^{12}$ is originally
obtained in \cite[Ex.6.8]{Re}. Indeed, we note that the equation
$y_{1}y_{2}=\left(\begin{array}{ccc}
x_{1} & x_{2} & x_{3}\end{array}\right)A^{\dagger}\left(\begin{array}{c}
x_{4}\\
x_{5}\\
x_{6}
\end{array}\right)$ given in \cite[p.24]{Re} corresponds to 
\[
st=-{\empty}^{t}\!xQ^{\dagger}y,
\]
where $A^{\dagger}$ and $Q^{\dagger}$ are the adjoint matrices of
$A$ and $Q$, respectively, and $A,y_{1},y_{2}$, $\left(\begin{array}{c}
x_{1}\\
x_{2}\\
x_{3}
\end{array}\right),\left(\begin{array}{c}
x_{4}\\
x_{5}\\
x_{6}
\end{array}\right)$ in ibid. correspond to $Q,-t,s,x,y$ here, respectively. We can verify
that 
\begin{align*}
\beta(x,y)= & -3{\empty}^{t}\!xQ^{\dagger}y\\
 & -2(q_{12}q_{21}+q_{13}q_{31}+q_{23}q_{32}+q_{22}^{2}+q_{22}q_{33}+q_{33}^{2}){\empty}^{t}\!xy.
\end{align*}
Moreover, we can also verify that all the other equations of $\mathfrak{Z}_{\mA}^{12}$
can be identified with the eight $4\times4$ Pfaffians in \cite[p.24, 25]{Re}.

\vspace{5pt}

\begin{lem}
\label{lem:The-affine-coordinateZ UFD}The following assertions hold:
\begin{enumerate}[$(1)$]

\item The scheme $\mathfrak{Z}_{\mA}^{12}$ is Gorenstein of codimension
$4$, and is normal.

\item The affine coordinate ring $R_{\mathfrak{Z}}$ of $\mathfrak{Z}_{\mA}^{12}$
is a UFD. 

\end{enumerate}
\end{lem}

\begin{proof}
(1). We can check that the $x_{i}$-, $y_{i}$-, $s$-, and $t$-charts
($i=1,2,3)$ of $\mathfrak{Z}_{\mA}^{12}$ are smooth and the union
of them is irreducible and is of codimension $4$ in the ambient affine
space. Using \cite{DGPS}, we can show that the ideal of $\mathfrak{Z}_{\mA}^{12}$
has the $9\times16$ graded minimal free resolution of length $4$.
Therefore, by the proofs of \cite[Prop.4.8]{Tak4}, we obtain (1).

(2). We can check that $(R_{\mathfrak{Z}})_{x_{1}}$ is a localization
of a polynomial ring. We can also check that $\mathfrak{Z}_{\mA}^{12}\cap\{x_{1}=0\}$
is normal. Indeed, on the $x_{2}$-, $x_{3}$-, $y_{2}$- and $y_{3}$-charts,
$\mathfrak{Z}_{\mA}^{12}\cap\{x_{1}=0\}$ is smooth, and $\mathfrak{Z}_{\mA}^{12}\cap\{x_{1}=x_{2}=x_{3}=y_{2}=y_{3}=0\}$
is of codimension $2$ in $\mathfrak{Z}_{\mA}^{12}\cap\{x_{1}=0\}$.
Therefore The proofs of \cite[Prop.4.9]{Tak4} work verbatim for (2).
\end{proof}
The main purpose of this section is to interpret $\mathfrak{Z}_{\mA}^{12}$
in the context of \cite{Tak3} as follows:
\begin{prop}
\label{prop:We-define-Sigma11}We define $\mathfrak{Z}_{\mP}^{11}$
to be the projective variety obtained from $\mathfrak{Z}_{\mA}^{12}$
by setting $w(s)=w(t)=2$ and all the other weights on the coordinates
as $1$ ( Note that $\mathfrak{Z}_{\mP}^{11}$ is contained in $\mP(1^{14},2^{2})$).
The variety $\mathfrak{Z}_{\mP}^{11}$ is isomorphic to the $\mQ$-Fano
variety $\Sigma$ as in \cite[Thm.1.1]{Tak3} associated to prime
$\mQ$-Fano $3$-folds of No.1.1 in \cite{Tak1}.
\end{prop}

\begin{proof}
Note that the projective variety $\overline{\mathfrak{Z}}:=\left\{ {\empty}^{t}\!xPy=0,{\empty}^{t}\!xQy=0\right\} \subset\mP^{13}$
with $P=E$ and ${\rm Tr\,}Q=0$ can be identified with $\overline{\Sigma}$
as in \cite[Subsec.4.1.1]{Tak3} with the obvious correspondence between
their coordinates. 

By the construction of $\Sigma$ from $\overline{\Sigma}$ as in \cite[Sec.5]{Tak3},
$\Sigma$ is isomorphic to $\overline{\Sigma}$ in codimension $1$
except the image $\overline{\Gamma}$ in $\overline{\Sigma}$ of the
exceptional divisor for the blow-up of two $1/2(1^{11})$-singularities. 

We show that a similar fact holds for $\mathfrak{Z}_{\mP}^{11}$ and
$\overline{\mathfrak{Z}}$. We consider the rational map $\pi\colon\mathfrak{Z}_{\mP}^{11}\dashrightarrow\mathfrak{\overline{Z}}$
which is the restriction of the projection from the $(s:t)$-line.
Note that $\pi$ is defined on $U:=\mathfrak{Z}_{\mP}^{11}\setminus\{\text{the\,}s\text{-, \ensuremath{t}-points\}. }$
We set $\overline{\Delta}:=\left\{ x=0,y^{\sharp}=0\right\} \cup\left\{ y=0,x^{\sharp}=0\right\} \subset\overline{\mathfrak{Z}}.$
From the equation of $\mathfrak{Z}_{\mP}^{11}$, we see that the closure
of the inverse image of $\overline{\Delta}|_{\pi(U)}$ by $\pi|_{U}$
is $\Delta:=\left\{ x=0,y^{\sharp}=0,t=0\right\} \cup\left\{ y=0,x^{\sharp}=0,s=0\right\} .$
By the equations of $\mathfrak{Z}_{\mP}^{11}$ and $\overline{\mathfrak{Z}}$,
we see that, if $x\not=0$ and $y\not=0$ on $\overline{\mathfrak{Z}}$,
then both $s$ and $t$ are recovered by the equations of $\mathfrak{Z}_{\mP}^{11}$.
Moreover, we see that points of $\{x=0\}\cup\{y=0$\} outside $\overline{\Delta}$
are not the $\pi$-images of points of $U$, and the fiber of $\pi|_{U}$
over a point ${\rm p}\in\overline{\Delta}$ is $1$-dimensional if
${\rm p}\not\in\{x=y=0\}$, and $2$-dimensional if ${\rm p}\in\{x=y=0\}$.
Therefore the $\pi|_{U}$-exceptional locus is $\Delta|_{U}$. By
the equation of $\mathfrak{Z}_{\mP}^{11}$ and \cite[Prop.4.4]{Tak3},
the locus $\overline{\Delta}$ coincides with ${\rm Sing}\,$$\mathfrak{\overline{Z}}$,
which is of codimension 3 in $\overline{\mathfrak{Z}}$. Therefore
the exceptional locus of $\pi|_{U}$ is of codimension 2 in $\mathfrak{Z}_{\mP}^{11}$,
which implies $\mathfrak{Z}_{\mP}^{11}$ is isomorphic to $\overline{\mathfrak{Z}}$
in codimension $1$ except $\{x=0\}\cup\{y=0$\}. Note that $\{x=0\}\cup\{y=0$\}
is identified with $\overline{\Gamma}$. Therefore, we see that $\Sigma$
and $\mathfrak{Z}_{\mP}^{11}$ are isomorphic in codimension $1$
identifying $\overline{\Sigma}$ and $\overline{\mathfrak{Z}}$. By
\cite[Thm.1.1]{Tak3}, $\Sigma$ is a $\mQ$-factorial $\mQ$-Fano
variety with Picard numbers 1, and, by Lemma \ref{lem:The-affine-coordinateZ UFD}
and the proof of \cite[Cor.4.10]{Tak4}, so is $\mathfrak{Z}_{\mP}^{11}$.
Therefore, $\Sigma$ and $\mathfrak{Z}_{\mP}^{11}$ are actually isomorphic
by \cite[Lem.5.5]{Tak2} as desired. 
\end{proof}
By Proposition \ref{prop:We-define-Sigma11} and \cite[Thm.1.1]{Tak3},
we obtain the following result of Gushel-Mukai type:
\begin{cor}
Any prime $\mQ$-Fano $3$-fold of No.1.1 in \cite{Tak1} is a weighted
complete intersection of $\mathfrak{Z}_{\mP}^{11}$ as in Proposition
\ref{prop:We-define-Sigma11} with respect to hypersurfaces of weight
$1$.
\end{cor}

Finally, we describe an ${\rm SL}_{3}$-action on $\mathfrak{Z}_{\mA}^{12}$
as follows:
\begin{prop}
\label{prop:Z12Gr} Let $\mA_{\mathfrak{Z}}^{16}$ be the affine subspace
$\{{\rm Tr}\,Q=0\}$ in the affine space whose coordinates are $s,\,t$
and the entries of $x,y,Q$. The affine space $\mathbb{A}_{\mathfrak{Z}}^{16}$
has the ${\rm SL}_{3}$-action preserving $\mathfrak{Z}_{\mA}^{12}$
defined by

\[
x\mapsto{\empty}^{t}\!gx,\,y\mapsto g^{-1}y,\,Q\mapsto g^{-1}Qg,s\mapsto s,t\mapsto t
\]
for $g\in{\rm SL}_{3}$.
\end{prop}

\section{\textbf{The $G_{2}^{(4)}$-Cluster variety as a specializations of
$\mathfrak{U}_{\mathbb{A}}^{14}$ \label{sec:The--Cluster-variety}}}

In this section, we clarify the relationship between the $G_{2}^{(4)}$-cluster
variety constructed in \cite{CD}, and the affine varieties $\mathfrak{F}_{\mathbb{A}}^{14}$
and $\mathfrak{Z}_{\mA}^{12}$. We review in our context the prime
$\mQ$-Fano 3-folds constructed in ibid.~from the $G_{2}^{(4)}$-cluster
variety. 

\subsection{The $G_{2}^{(4)}$-cluster variety $\mathfrak{Cl}_{\mA}^{10}$}

In \cite{CD}, the $G_{2}^{(4)}$-cluster variety is defined in the
affine $16$-space with coordinates 
\begin{align*}
 & \theta_{i}\,(1\leq i\leq4),\,\theta_{23},\,\theta_{41},\\
 & A_{j}\,(1\leq j\leq4),A_{kl}\,((k,l)=(1,2),\,(2,3),\,(3,4),\,(4,1)),\\
 & \lambda_{13},\,\lambda_{24}.
\end{align*}
In this paper, however, we call a smaller variety as the $G_{2}^{(4)}$-cluster
variety. By the big table \cite{bigtables}, we observe that $A_{12}$
and $A_{34}$ are always nonzero constants when a $\mQ$-Fano 3-fold
is constructed from a weighted projectivization of the $G_{2}^{(4)}$-cluster
variety. Then, replacing $A_{1}$, $A_{2}$, $A_{3}$, $A_{4}$, $\lambda_{13}$,
$\lambda_{24}$ with $A_{12}^{-1}A_{1}$, $A_{12}^{-1}A_{2}$, $A_{34}^{-1}A_{3}$,
$A_{34}^{-1}A_{4}$, $A_{12}^{-1}A_{34}^{-1}\lambda_{13}$, $A_{12}^{-1}A_{34}^{-1}\lambda_{24}$,
we see that it is possible to set $A_{12}=A_{34}=1$. 
\begin{defn}
\textit{\label{def:G24CL}The $G_{2}^{(4)}$-cluster variety} $\mathfrak{Cl}_{\mA}^{10}$
is defined as a subvariety of the affine $14$-space with coordinates
\begin{align*}
 & \theta_{i}\,(1\leq i\leq4),\theta_{23},\,\theta_{41},\\
 & A_{j}\,(1\leq j\leq4),A_{kl}\,((k,l)=(2,3),(4,1)),\\
 & \lambda_{13},\,\lambda_{24}
\end{align*}
by setting for the equations of $\mathfrak{F}_{\mathbb{A}}^{22}$
as follows:
\begin{align*}
 & x={\empty}^{t}\!\left(\begin{array}{ccc}
\theta_{4} & \theta_{1} & A_{23}\end{array}\right),y={\empty}^{t}\!\left(\begin{array}{ccc}
A_{41} & \theta_{2} & \theta_{3}\end{array}\right),\\
 & P=\left(\begin{array}{ccc}
-A_{4} & 0 & 0\\
0 & 0 & 1\\
-\lambda_{24} & -A_{2} & 0
\end{array}\right),Q=\left(\begin{array}{ccc}
0 & 1 & 0\\
-A_{1} & 0 & 0\\
-\lambda_{13} & 0 & -A_{3}
\end{array}\right),\\
 & s=-\theta_{23},t=-\theta_{41}.
\end{align*}
\end{defn}

We can immediately check that this definition of the $G_{2}^{(4)}$-cluster
variety coincides with that in \cite[Subsec.1.2.2]{CD} when $A_{12}=A_{34}=1$. 

By an elementary calculation, we have the following:
\begin{prop}
\label{prop:All-the-equations}All the equations of $\mathfrak{Cl}_{\mA}^{10}$
are weighted homogeneous if and only if the following conditions for
the weights on the coordinates hold:
\begin{align*}
 & w(A_{1})=-2w(\theta_{1})+w(\theta_{2})+w(\theta_{41}),\,w(A_{2})=w(\theta_{1})-2w(\theta_{2})+w(\theta_{23}),\\
 & w(A_{3})=-2w(\theta_{3})+w(\theta_{4})+w(\theta_{23}),\,w(A_{4})=w(\theta_{3})-2w(\theta_{4})+w(\theta_{41}),\\
 & w(A_{23})=w(\theta_{2})+w(\theta_{3})-w(\theta_{23}),\,w(A_{41})=w(\theta_{1})+w(\theta_{4})-w(\theta_{41}),\\
 & w(\lambda_{13})=-w(\theta_{1})-w(\theta_{3})+w(\theta_{23})+w(\theta_{41}),\\
 & w(\lambda_{24})=-w(\theta_{2})-w(\theta_{4})+w(\theta_{23})+w(\theta_{41}).
\end{align*}
\end{prop}

\vspace{5pt}

In the following subsections, we only consider the weights on the
coordinates as in Proposition \ref{prop:All-the-equations} such that
the weighted projectivizations of $\mathfrak{Cl}_{\mA}^{10}$ itself
or its subvarieties associated to these produce prime $\mQ$-Fano
3-folds. The list of such weights are presented in \cite{bigtables}. 

\subsection{The maximal case\label{subsec:The-maximal-case}}

By the big table \cite{bigtables}, we observe that all the weights
on the coordinates of the $G_{2}^{(4)}$-cluster variety $\mathfrak{Cl}_{\mA}^{10}$
are positive only in the two cases No.5530 and No.11455 of \cite{GRDB}.
We may verify the following by a straightforward calculation.
\begin{prop}
\label{prop:We-define-theCL}The equations of the affine variety $\mathfrak{Cl}_{\mA}^{10}$
are presented in the format of the equations of $\mathfrak{U}_{\mA}^{14}$
by setting 

\begin{align*}
 & x={\empty}^{t}\!\left(\begin{array}{ccc}
\theta_{1}-A_{2}A_{23} & \theta_{4}-A_{3}A_{23} & A_{23}\end{array}\right),y={\empty}^{t}\!\left(\begin{array}{ccc}
A_{41} & 1/2(\theta_{3}-\theta_{2}) & 1/2(\theta_{3}+\theta_{2})\end{array}\right),\\
 & P=\left(\begin{array}{ccc}
0 & 1 & 1\\
-A_{4} & 0 & 0\\
-\lambda_{24}-A_{3}A_{4} & 2A_{2} & 0
\end{array}\right),Q=\left(\begin{array}{ccc}
-A_{1} & 0 & 0\\
0 & -1 & 1\\
-\lambda_{13}-A_{1}A_{2} & -2A_{3} & 0
\end{array}\right),\\
 & s=\theta_{23},t=2\theta_{41}.
\end{align*}
\end{prop}

By Proposition \ref{prop:All-the-equations} and the weights on the
coordinates of $\mathfrak{Cl}_{\mA}^{10}$ for No.5530 and No.11455,
we see that all the entries of $x,y,P,Q$ as in Proposition \ref{prop:We-define-theCL}
are weighted homogeneous. Therefore we have the following:
\begin{cor}
The affine variety $\mathfrak{Cl}_{\mA}^{10}$ is isomorphic to the
subvariety of $\mathfrak{U}_{\mathbb{A}}^{14}$ 
\begin{equation}
\mathfrak{U}_{\mathbb{A}}^{14}\cap\{a_{11}=0,b_{11}=1,b_{12}=0,b_{21}=0\}.\label{eq:max}
\end{equation}
Moreover, the isomorphism is weighted homogeneous for No.5530 and
No.11455 with the corresponding weights on the coordinates of \textup{(\ref{eq:max}).}
\end{cor}

\subsection{The case where only one coordinate is a nonzero constant\label{subsec:The-case-onlyone}}

By \cite{bigtables}, there are 12 classes of prime $\mQ$-Fano 3-folds
for which only one coordinate of $\mathfrak{Cl}_{\mA}^{10}$ (actually
$A_{1}$,$A_{3}$ or $A_{4}$) is a nonzero constant. By symmetry,
we may identify the three cases, so we have only to consider the case
that $A_{4}$ is a nonzero constant; this is the case for No.1169,
1182, 5860, 5870, 16228. We denote by $\mathfrak{Cl}_{\mA,A_{4}}^{9}$
the subvariety of $\mathfrak{Cl}_{\mA}^{10}$ with $A_{4}=-1$. By
suitable changes of scales of coordinates, we may assume that $A_{4}=-1$. 

We may verify the following by a straightforward calculation.
\begin{prop}
\label{prop:We-define-theCL9}The equations of the affine variety
$\mathfrak{Cl}_{\mA,A_{4}}^{9}$ are presented in the format of the
equations of $\mathfrak{U}_{\mA}^{14}$ by setting
\begin{align*}
 & x={\empty}^{t}\!\left(\begin{array}{ccc}
A_{41} & \theta_{2} & \theta_{3}\end{array}\right),y={\empty}^{t}\!\left(\begin{array}{ccc}
A_{23} & \theta_{1} & \theta_{4}-1/2\lambda_{24}A_{23}\end{array}\right),\\
 & P=\left(\begin{array}{ccc}
-1/2\lambda_{24} & 0 & 1\\
-A_{2} & 0 & 0\\
0 & 1 & 0
\end{array}\right),Q=\left(\begin{array}{ccc}
-\lambda_{13} & -A_{1} & 0\\
1/2\lambda_{24} & 0 & 1\\
-A_{3} & 0 & 0
\end{array}\right),\\
 & s=-\theta_{41},t=\theta_{23}.
\end{align*}
\end{prop}

By Proposition \ref{prop:All-the-equations} and the assumption that
$A_{4}=-1$, we see that the entry $\theta_{4}-1/2\lambda_{24}A_{23}$
of $y$ is weighted homogeneous. Therefore we have the following:
\begin{cor}
The affine variety $\mathfrak{Cl}_{\mA,A_{4}}^{9}$ is isomorphic
to the subvariety of $\mathfrak{U}_{\mathbb{A}}^{14}$ 
\begin{equation}
\mathfrak{U}_{\mathbb{A}}^{14}\cap\{b_{11}=b_{12}=c_{11}=c_{22}=0,c_{12}=1\}.\label{eq:semimax}
\end{equation}
Moreover, the isomorphism is weighted homogeneous with the corresponding
weights on the  coordinates of \textup{(\ref{eq:semimax}).}
\end{cor}

\vspace{10pt}

In the following subsections, we consider the cases in which exactly
two of the coordinates are nonzero constants. There are 4 cases; the
set of the two nonzero constant coordinates is $\{A_{3},A_{4}\}$
, $\{A_{1},A_{4}\}$, $\{A_{2},A_{3}\}$, or $\{A_{1},A_{3}\}$. For
them, we may assume that the nonzero constants are $-1$. We denote
by $\mathfrak{Cl}_{\mA,A_{i}\!A_{j}}^{8}$ the subvariety of $\mathfrak{Cl}_{\mA}^{10}$
with $A_{i}=A_{j}=-1$. We will treat the 4 cases separately in the
sequel. 

\subsection{The case where $A_{3}$ and $A_{4}$ are nonzero constants\newline--another
appearance of $\mathfrak{S}_{\mA}^{8}$--\label{subsec:The-case-A3A4}}

This is the case for 12 classes in \cite{GRDB}. In this subsection,
we see another appearance of the affine variety $\mathfrak{S}_{\mA}^{8}$.

We may verify the following by a straightforward calculation.
\begin{prop}
\label{prop:We-define-theCLS8}The equations of the affine variety
$\mathfrak{Cl}_{\mA,A_{3}\!A_{4}}^{8}$ are presented in the format
of the equations of $\mathfrak{S}_{\mA}^{8}$ by setting
\begin{align*}
 & x={\empty}^{t}\!\left(\begin{array}{ccc}
A_{41} & \theta_{2} & \theta_{3}-2/3\lambda_{13}A_{41}+1/3\lambda_{24}\theta_{2}\end{array}\right),\\
 & y={\empty}^{t}\!\left(\begin{array}{ccc}
\theta_{1} & A_{23} & \theta_{4}+1/3\lambda_{13}\theta_{1}-2/3\lambda_{24}A_{23}\end{array}\right),\\
 & P=\left(\begin{array}{ccc}
1/3\lambda_{13} & -1/3\lambda_{24} & 1\\
-1/3\lambda_{24} & -A_{2} & 0\\
1 & 0 & 0
\end{array}\right),Q=\left(\begin{array}{ccc}
-A_{1} & -1/3\lambda_{13} & 0\\
-1/3\lambda_{13} & 1/3\lambda_{24} & 1\\
0 & 1 & 0
\end{array}\right),\\
 & s=-\theta_{41},t=-\theta_{23}.
\end{align*}
\end{prop}

By Proposition \ref{prop:All-the-equations} and the assumption that
$A_{3}=A_{4}=-1$, we see that the entries $\theta_{3}-2/3\lambda_{13}A_{41}+1/3\lambda_{24}\theta_{2}$
and $\theta_{4}+1/3\lambda_{13}\theta_{1}-2/3\lambda_{24}A_{23}$
of $x$ and $y$ are weighted homogeneous. Therefore we have the following:
\begin{cor}
The affine variety $\mathfrak{Cl}_{\mA,A_{3}\!A_{4}}^{8}$ is isomorphic
to $\mathfrak{S}_{\mA}^{8}$. Moreover, the isomorphism is weighted
homogeneous with the corresponding weights on the coordinates of the
equations as in Proposition \ref{prop:NewS8}.
\end{cor}

\begin{example}
\label{exa:We-consider-theEx20652}We consider the weighted projectivization
$\mathfrak{S}_{\mP}^{7}$ of $\mathfrak{S}_{\mathbb{A}}^{8}$ in $\mP(1^{8},2^{4})$
by putting the weights on the coordinates as follows:
\[
w(u_{i})=1\,(1\leq i\leq4),w(v_{j})=2\,(0\leq j\leq3),w(d_{k})=1\,(0\leq k\leq3).
\]
By \cite{CD}, this produces a prime $\mQ$-Fano 3-fold $X$ of No.20652
as its weighted complete intersection by one hypersurface of weight
2 and three hypersurfaces of weight 1. It is easy to verify that the
the set of three $1/2\,(1,1,1)$-singularities of $X$ is the intersection
between the twisted cubic curve $\gamma:=\{U=O,D=0,\wedge^{2}V=\bm{o}\}$
and $X$.
\end{example}

We refine Example \ref{exa:We-consider-theEx20652}. In the paper
\cite{Tak1}, we obtain two different classes No.4.1 and No.5.4 of
prime $\mQ$-Fano threefolds of No.20652. In the following proposition,
we determine which class a prime $\mQ$-Fano $3$-fold $X$ constructed
in Example \ref{exa:We-consider-theEx20652} belongs to:
\begin{prop}
A prime $\mQ$-Fano $3$-fold $X$ constructed in Example \ref{exa:We-consider-theEx20652}
belongs to the class No.4.1 of \cite{Tak1} for any $1/2\,(1,1,1)$-singularity
of $X$.
\end{prop}

\begin{proof}
We choose a $1/2\,(1,1,1)$-singularity $\mathsf{p}$ and construct
the Sarkisov link starting from the blow-up at $\mathsf{p}$. By \cite[Tables 4 and 5]{Tak1},
we can distinguish between prime $\mQ$-Fano threefolds of No.4.1
and No.5.4 by the dimension of forms of weight 2 vanishing at $\mathsf{p}$
with weighted multiplicity $\geq$ $3$; the dimension is $3$ for
No.4.1, and is 4 for No.5.4, which can be verified by looking at the
rational map $X\dashrightarrow X'$ in the Sarkisov link (note that
this map is defined by the linear system $|-2K_{X}-3\mathsf{p}|)$.
By the group action on $\mathfrak{S}_{\mA}^{8}$ described in Subsection
\ref{subsec:-action}, we may assume that $\mathsf{p}$ is the $v_{0}$-point
(the point whose coordinates except $v_{0}$ are zero). Then, by the
equation (\ref{eq:S8}), we see that 4 forms of weight 2, $u_{1}^{2},u_{1}u_{3},u_{3}^{2},v_{3}$
are all such forms and they are linearly independent on $\mathfrak{S}_{\mA}^{8}$.
We have only to show that they are still linearly independent on $X$.
Assume the contrary. Then we may assume that the hypersurface $Q$
of weight 2 cutting $X$ from $\mathfrak{S}_{\mP}^{7}$ is $\{av_{3}+bu_{1}^{2}+cu_{1}u_{3}+du_{3}^{2}=0\}$
with some $a,b,c,d\in\mC$. Then $Q$ intersects the twisted cubic
curve $\gamma$ only at the $v_{0}$-point, or contains $\gamma$.
This is a contradiction since $Q$ must intersects $\gamma$ at three
distinct points for $X$ to have three $1/2\,(1,1,1)$-singularities.
\end{proof}
\begin{rem}
In \cite{Tak4}, we will construct a prime $\mQ$-Fano threefold of
No.5.4 via another key variety $\mathscr{H}_{\mA}^{13}$.
\end{rem}

\subsection{The case where $A_{1}$ and $A_{4}$, or $A_{2}$ and $A_{3}$ are
nonzero constants\newline --the affine variety $\mathfrak{T}_{\mA}^{8}$
with an ${\rm SL_{2}}$-action--\label{subsec:The-case-AA1A4}}

This is the case for 16 classes in \cite{GRDB}. It is easy to see
that we may identify $\mathfrak{Cl}_{\mA,A_{1}\!A_{4}}^{8}$ and $\mathfrak{Cl}_{\mA,A_{2}\!A_{3}}^{8}$
changing the suffixes of the coordinates. Hence we only consider $\mathfrak{Cl}_{\mA,A_{1}\!A_{4}}^{8}$
in this subsection. This is the case for 12 classes.

We may verify the following by a straightforward calculation.
\begin{prop}
\label{prop:We-define-theCLT8}The equations of the affine variety
$\mathfrak{Cl}_{\mA,A_{1}\!A_{4}}^{8}$ are presented in the format
of the equations of $\mathfrak{U}_{\mA}^{14}$ by the setting as in
Proposition \ref{prop:We-define-theCL9} with $A_{1}=-1$. 
\end{prop}

Now we arrive at a new interpretation of $\mathfrak{Cl}_{\mA,A_{1}\!A_{4}}^{8}$;
we set
\begin{align*}
 & w_{1}:=\theta_{1}-1/3\lambda_{13}A_{23},\,w_{2}:=-\theta_{4}+1/3\lambda_{24}A_{23},\\
 & z_{1}:=A_{41},\,z_{2}:=\theta_{2},\,z_{3}:=\theta_{3},\\
 & f_{0}:=-A_{3},\,f_{1}:=1/3\lambda_{13},\,f_{2}:=1/3\lambda_{24},f_{3}:=-A_{2},\\
 & s:=-\theta_{41},t:=\theta_{23},u:=A_{23}.
\end{align*}

By Proposition \ref{prop:All-the-equations} and the assumption that
$A_{1}=A_{4}=-1$, we see that the entries $\theta_{1}-1/3\lambda_{13}A_{23}$
and $-\theta_{4}+1/3\lambda_{24}A_{23}$ are weighted homogeneous.
Moreover, we set

\begin{align*}
 & w:=\left(\begin{array}{c}
w_{1}\\
w_{2}
\end{array}\right),\,Z:=\left(\begin{array}{cc}
z_{1} & -z_{2}\\
z_{3} & -z_{1}
\end{array}\right),\,z:=\left(\begin{array}{c}
z_{2}\\
-2z_{1}\\
z_{3}
\end{array}\right),\\
 & F:=\left(\begin{array}{ccc}
f_{2} & f_{1} & f_{0}\\
f_{3} & f_{2} & f_{1}
\end{array}\right),\,F^{\dagger}:=\left(\begin{array}{cc}
-f_{1} & f_{0}\\
2f_{2} & -2f_{1}\\
-f_{3} & f_{2}
\end{array}\right).
\end{align*}
 
\begin{prop}
\label{prop:The-affine-varietyT8}The affine variety $\mathfrak{Cl}_{\mA,A_{1}\!A_{4}}^{8}$
is isomorphic to the affine subvariety $\mathfrak{T}_{\mA}^{8}$ in
the affine $8$-space $\mA_{\mathfrak{T}}^{12}$ with the coordinates
\[
w_{1},w_{2},z_{1},z_{2},z_{3},s,t,u,f_{0},f_{1},f_{2},f_{3}
\]
 defined by the following equations: 
\begin{align*}
 & Zw+uFz=\bm{o},\\
 & tw=ZFz,\\
 & tu=\det Z,\\
 & sz=-2u^{2}\wedge^{2}({\empty}^{t}\!F)+uF^{\dagger}w+\left(\begin{array}{c}
w_{1}^{2}\\
-2w_{1}w_{2}\\
w_{2}^{2}
\end{array}\right),\\
 & st=-1/2\,u(\wedge^{2}F^{\dagger})z+\left(\begin{array}{cc}
w_{2} & -w_{1}\end{array}\right)Fz.
\end{align*}
Moreover, the isomorphism is weighted homogeneous with the corresponding
weights on the coordinates of these equations.
\end{prop}

Though the equations of the affine variety $\mathfrak{T}_{\mA}^{8}$
look complicated, it turn out to be suitable to see an ${\rm SL}_{2}$-action
on $\mathfrak{T}_{\mA}^{8}$ as follows:
\begin{prop}
\label{prop:T8Orbit}The affine space $\mathbb{A}_{\mathfrak{T}}^{12}$
has the ${\rm SL}_{2}$-action preserving $\mathfrak{T}_{\mA}^{8}$
defined by

\[
Z\mapsto gZg^{-1},w\mapsto gw,F\mapsto gF\widehat{g}^{-1},s\mapsto s,\,t\mapsto t,\,u\mapsto u,
\]
where $g\in{\rm SL}_{2}$, and the definition of $\widehat{g}$ for
$g\in{\rm SL}_{2}$ is as in Proposition \ref{prop:S8Orbit}.
\end{prop}

\begin{proof}
We only note that $z$ is mapped to $\widehat{g}z$ by the ${\rm SL}_{2}$-action. 
\end{proof}

\subsection{The case where $A_{1}$ and $A_{3}$ are nonzero constants-- a subvariety
of $\mathfrak{Z}_{\mathbb{A}}^{12}$--\label{subsec:The-case-A1A3}}

We may verify the following by a straightforward calculation.
\begin{prop}
\label{prop:We-defineZ12 sub}The equations of the affine variety
$\mathfrak{Cl}_{\mA}^{10}$ are presented in the format of the equations
of $\mathfrak{Z}_{\mA}^{12}$ by setting

\begin{align*}
 & x={\empty}^{t}\!\left(\begin{array}{ccc}
\theta_{4} & \theta_{1} & A_{23}\end{array}\right),y={\empty}^{t}\!\left(\begin{array}{ccc}
\theta_{2} & A_{41} & \theta_{3}-\lambda_{13}A_{41}\end{array}\right),\\
 & P=E,\,Q=\left(\begin{array}{ccc}
-1/3\lambda_{13} & -A_{4} & 0\\
0 & 2/3\lambda_{13} & 1\\
-A_{2} & -\lambda_{24} & -1/3\lambda_{13}
\end{array}\right),\\
 & s=\theta_{23},t=-\theta_{41}.
\end{align*}
\end{prop}

By Proposition \ref{prop:All-the-equations} and the assumption that
$A_{1}=A_{3}=-1$, we see that the entry $\theta_{3}-\lambda_{13}A_{41}$
of $y$ is weighted homogeneous. Therefore we have the following:
\begin{cor}
The affine variety $\mathfrak{Cl}_{\mA,A_{1}\!A_{3}}^{8}$ is isomorphic
to the subvariety 
\[
\{q_{22}=-2q_{11},q_{33}=q_{11},q_{13}=q_{21}=0,q_{23}=1\}
\]
 of $\mathfrak{Z}_{\mathbb{A}}^{12}$. Moreover, the isomorphism is
weighted homogeneous with the corresponding weights on the coordinates
of the equations as in Definition \ref{def:We-set-Z12}.
\end{cor}

This is the case for 16 classes in \cite{GRDB}.

\subsection{The case with the most number of nonzero constant coordinates \newline
--the affine variety $\mathfrak{B}_{\mA}^{6}$ with an ${\rm SL}_{2}$-action--
\label{subsec:The-case-withtwoparam.}}

By \cite{bigtables}, the number of nonzero constant coordinates is
at most 3, and, if the number is 3, then the nonzero constant coordinates
are always $A_{1}$, $A_{3}$ and $A_{4}$. This is the case for 35
classes in \cite{GRDB}. In this case, we may assume that $A_{1}=A_{3}=A_{4}=-1$.
Then the subvariety of $\mathfrak{Cl}_{\mA}^{10}$ with $A_{1}=A_{3}=A_{4}=-1$
is a subvariety of $\mathfrak{Cl}_{\mA,A_{1}\!A_{3}}^{8}$, $\mathfrak{Cl}_{\mA,A_{1}\!A_{4}}^{8}$
and $\mathfrak{Cl}_{\mA,A_{3}\!A_{4}}^{8}$.

Here we consider that this is a subvariety of $\mathfrak{Cl}_{\mA,A_{3}\!A_{4}}^{8}$
and then identifying $\mathfrak{Cl}_{\mA,A_{3}\!A_{4}}^{8}$ with
$\mathfrak{S}_{\mA}^{8}$, we describe $\mathfrak{Cl}_{\mA,A_{3}\!A_{4}}^{8}$
as the subvariety $\mathfrak{S}_{\mA}^{8}\cap\{d_{3}=-1\}$ of $\mathfrak{S}_{\mA}^{8}$.
Moreover, defining the new coordinates $D_{0},D_{1},U_{1},U_{2}$
corresponding to $d_{0},d_{1},u_{1},u_{2}$ by 

\[
D_{0}=d_{0}+3d_{1}d_{2}+2d_{2}^{3},D_{1}=d_{1}+d_{2}^{2},U_{1}=u_{1}-d_{2}u_{3},U_{2}=u_{2}-d_{2}u_{4},
\]
we see by a straightforward calculation that $\mathfrak{S}_{\mA}^{8}\cap\{d_{3}=-1\}$
is isomorphic to the cone over $\mathfrak{S}_{\mA}^{8}\cap\{d_{2}=0,\,d_{3}=-1\}.$
The isomorphism is weighted homogeneous with the corresponding weights
on the coordinates since $D_{0},D_{1},U_{1},U_{2}$ are weighted homogeneous
if $A_{1}=A_{3}=A_{4}=-1$ by Proposition \ref{prop:All-the-equations}.

We set 
\[
\mathfrak{B}_{\mA}^{6}:=\mathfrak{S}_{\mA}^{8}\cap\{d_{2}=0,\,d_{3}=-1\}.
\]

Let $\rho_{\mathfrak{B}}:\widehat{\mathfrak{B}}\to\{d_{2}=0,d_{3}=-1\}$
be the base change of $\rho_{\mathfrak{S}}$ as in Subsection \ref{subsec:-fibration-8dim}
by the inclusion map $\{d_{2}=0,d_{3}=-1\}\hookrightarrow\mA_{D}^{4}$.
By Proposition \ref{prop:P1P1P1FibS8}, we have immediately the following
noting that the cone over the tangential scroll (\ref{eq:tangential})
of the twisted cubic restricts to the affine cuspidal cubic curve
$\{d_{2}=0,d_{3}=-1,d_{0}^{2}=4d_{1}^{3}\}$:
\begin{prop}
Let $\mathsf{p}$ be a point of $\{d_{2}=0,d_{3}=-1\}$ and $F_{\mathsf{p}}$
the $\rho_{\mathfrak{B}}$-fiber over $\mathsf{p}.$ We identify $\{d_{2}=0,d_{3}=-1\}$
with the affine $2$-space with $d_{0},d_{1}$ as the coordinates.

If $\mathsf{p}=0$, then $F_{\mathsf{p}}\simeq\mP^{1,1,1}.$

If $\mathsf{p}$ belongs to $\{d_{0}^{2}=4d_{1}^{3}\}\setminus\{0\}$
, then $F_{\mathsf{p}}\simeq\mP^{1}\times\mathtt{Q}$.

If $\mathsf{p}$ does not belong to $\{d_{0}^{2}=4d_{1}^{3}\}$, then
$F_{\mathsf{p}}\simeq\mP^{1}\times\mP^{1}\times\mP^{1}.$
\end{prop}

We can also easily check the following:
\begin{prop}
The variety $\mathfrak{B}_{\mA}^{6}$ has an ${\rm SL_{2}}$-action
which is the restriction of the ${\rm SL_{2}^{II}}$-action of $\mathfrak{S}_{\mA}^{8}$.
\end{prop}

\subsection{Summary of the results}

In the following table, we summarize the numbers of classes of prime
$\mQ$-Fano 3-folds which are obtained from the $G_{2}^{(4)}$-cluster
variety $\mathfrak{Cl}_{\mA}^{10}$ and its subvarieties:

\begin{table}[H]
\begin{center}

\begin{tabular}{|c|c|c|}
\hline 
Subsection & constant coordinates & $\sharp$ of classes\tabularnewline
\hline 
\hline 
\ref{subsec:The-maximal-case} & None & $2$\tabularnewline
\hline 
\ref{subsec:The-case-onlyone} & $A_{1}$, $A_{3}$, or $A_{4}$ & $12$\tabularnewline
\hline 
\ref{subsec:The-case-A3A4} & $A_{3}$ and $A_{4}$ & $12$\tabularnewline
\hline 
\multirow{2}{*}{\ref{subsec:The-case-AA1A4}} & $A_{1}$ and $A_{4}$, or & \multirow{2}{*}{$16$}\tabularnewline
 & $A_{2}$ and $A_{3}$~~~~~~ & \tabularnewline
\hline 
\ref{subsec:The-case-A1A3} & $A_{1}$ and $A_{3}$ & $16$\tabularnewline
\hline 
\ref{subsec:The-case-withtwoparam.} & $A_{1}$, $A_{3}$ and $A_{4}$ & $35$\tabularnewline
\hline 
\end{tabular}

\end{center}\caption{prime $\mQ$-Fano 3-fold obtained from $\mathfrak{Cl}_{\mA}^{10}$}
\end{table}

\section{\textbf{Affine variety }$\mathfrak{P}_{\mA}^{23}$\textbf{ defined
by type II$_{2}$ unprojection due to S.~Papadakis\label{sec:Affine-variety-Papadakis}}}

In the paper \cite{P3}, Papadakis constructs affine varieties via
type ${\rm II}_{2}$ unprojection and proves they are Gorenstein (see
\cite[Thm.2.15]{P3}). Moreover, in \cite{P4}, he provides explicit
descriptions of parts of their equations, and a full description of
the equation in one particular case (see \cite[Sec.4]{P4}), for which
we denote by $\mathfrak{P}_{\mA}^{23}$ the affine variety he constructs.
The affine variety $\mathfrak{P}_{\mA}^{23}$ is a subvariety of codimension
4 in the affine $27$-space with the coordinates implemented in the
following format:

\begin{align*}
 & A^{k}=\left(\begin{array}{ccc}
0 & a_{12}^{k} & a_{13}^{k}\\
-a_{12}^{k} & 0 & a_{23}^{k}\\
-a_{13}^{k} & -a_{23}^{k} & 0
\end{array}\right),\,B^{k}=\left(\begin{array}{ccc}
b_{11}^{k} & b_{12}^{k} & b_{13}^{k}\\
b_{12}^{k} & b_{22}^{k} & b_{23}^{k}\\
b_{13}^{k} & b_{23}^{k} & b_{33}^{k}
\end{array}\right)\,(k=1,2),\\
 & {\empty}^{t}\!X=\left(\begin{array}{ccc}
X_{1} & X_{2} & X_{3}\end{array}\right),{\empty}^{t}\!Y=\left(\begin{array}{ccc}
Y_{1} & Y_{2} & Y_{3}\end{array}\right),\\
 & s_{0},\,s_{1},\,z
\end{align*}
(here we use a slightly different notation from Papadakis' one; we
denote by $X_{i}$, $Y_{i}$, $a_{ij}^{k}$ and $b_{ij}^{k}$ his
$x_{i}$, $y_{i}$, $A_{ij}^{k}$ and $B_{ij}^{k}$ respectively).
In the following proposition, we clarify that the variety $\mathfrak{P}_{\mA}^{23}$
is closely related with the variety $\mathfrak{F}_{\mA}^{22}$. The
assertion follows by a straightforward calculation. 

We set

\begin{align*}
 & D_{X}=\det\left(\begin{array}{ccc}
X_{1} & X_{2} & X_{3}\\
a_{23}^{1} & -a_{13}^{1} & a_{12}^{1}\\
a_{23}^{2} & -a_{13}^{2} & a_{12}^{2}
\end{array}\right),\,D_{Y}=\det\left(\begin{array}{ccc}
Y_{1} & Y_{2} & Y_{3}\\
a_{23}^{1} & -a_{13}^{1} & a_{12}^{1}\\
a_{23}^{2} & -a_{13}^{2} & a_{12}^{2}
\end{array}\right),\\
 & {\empty}^{t}\!\bm{a}_{1}=\left(\begin{array}{ccc}
-a_{23}^{1} & a_{13}^{1} & -a_{12}^{1}\end{array}\right),{\empty}^{t}\!\bm{a}_{2}=\left(\begin{array}{ccc}
-a_{23}^{2} & a_{13}^{2} & -a_{12}^{2}\end{array}\right).
\end{align*}

\begin{prop}
\label{prop:The-affine-varietyP23}The affine variety $\mathfrak{P}_{\mA}^{23}$
is transformed over the locus $\{z\not=0\}$ to the cone over $\mathfrak{F}_{\mA}^{22}$
by the following correspondence between coordinates: 
\begin{align*}
 & x=\sqrt{z}X+Y,\,y=-\sqrt{z}X+Y,\\
 & P=\frac{1}{2\sqrt{z}}A_{1}+B_{1},\,Q=\frac{1}{2\sqrt{z}}A_{2}+B_{2},\\
 & s=-\sqrt{z}s_{0}-s_{1}-\frac{3}{4\sqrt{z}}D_{X}+\frac{1}{4z}D_{Y}+\frac{1}{2\sqrt{z}}({\empty}^{t}\!\bm{a}_{1}B^{2}Y-{\empty}^{t}\!\bm{a}_{2}B^{1}Y)-\frac{1}{2}({\empty}^{t}\!\bm{a}_{1}B^{2}X-{\empty}^{t}\!\bm{a}_{2}B^{1}X),\\
 & t=\sqrt{z}s_{0}-s_{1}+\frac{3}{4\sqrt{z}}D_{X}+\frac{1}{4z}D_{Y}-\frac{1}{2\sqrt{z}}({\empty}^{t}\!\bm{a}_{1}B^{2}Y-{\empty}^{t}\!\bm{a}_{2}B^{1}Y)-\frac{1}{2}({\empty}^{t}\!\bm{a}_{1}B^{2}X-{\empty}^{t}\!\bm{a}_{2}B^{1}X),
\end{align*}
where the l.h.s.~and the r.h.s.~of the equalities correspond to
$\mathfrak{F}_{\mA}^{22}$ and $\mathfrak{P}_{\mA}^{23}$ respectively,
and $z$ is also the free coordinate of the cone over $\mathfrak{F}_{\mA}^{22}$.
\end{prop}

\begin{rem}
\begin{enumerate}

\item The definition of $b$ in \cite[p.2203, (4.2)]{P4} is slightly
incorrect; the coefficients of $B_{ij}^{1}{\rm ad}B_{ij}^{2}$ and
$B_{ij}^{2}{\rm ad}B_{ij}^{1}$ in the last part of r.h.s. of (4.2)
should be 2 when $i\not=j$. 

\item In \cite[Sec.5]{P4}, two candidates of prime $\mQ$-Fano 3-folds
of anti-canonical codimension 4 with type ${\rm II}_{1}$ projections
are constructed from $\mathfrak{P}_{\mA}^{23}$. In her phD thesis
submitted to Warwick University \cite{Tay}, Taylor constructs more
such candidates. The remaining problem is to show that their examples
have Picard number 1. 

\item By Proposition \ref{prop:The-affine-varietyP23} and Corollary
\ref{cor:P1P1P1}, the suitable partial projectivization of $\mathfrak{P}_{\mA}^{23}$
has a $\mP^{1}\times\mP^{1}\times\mP^{1}$-fibration.

\end{enumerate}
\end{rem}

\end{document}